\documentclass[a4paper,12pt, twoside, reqno]{amsart}
\usepackage{amssymb}
\usepackage{amsmath}
\usepackage{mathtools}
\usepackage{amsthm} 
\usepackage{tikz}
\usepackage{float} 

\usepackage{caption}
\DeclareCaptionFont{xipt}{\fontsize{10}{13}\mdseries}
\usepackage[font=xipt]{caption}

\usetikzlibrary{calc}
\usetikzlibrary{external}
\usepackage{xcolor}
\definecolor{bluegray}{rgb}{0.4, 0.6, 0.8}
\definecolor{blizzardblue}{rgb}{0.67, 0.9, 0.93}
\definecolor{myblue}{rgb}{0.137,0.466,0.741}
\definecolor{mygreen}{rgb}{0.137,0.466,0.1}
\definecolor{myred}{rgb}{0.837,0.166,0.141}

\definecolor{cerulean}{rgb}{0.0, 0.48, 0.65}
\definecolor{darkred}{rgb}{0.55, 0.0, 0.0}
\definecolor{tangelo}{rgb}{0.98, 0.3, 0.0}
\definecolor{plum}{rgb}{0.56, 0.27, 0.52}
\definecolor{olive}{rgb}{0.5, 0.5, 0.0}
\definecolor{applegreen}{rgb}{0.55, 0.71, 0.0}

\numberwithin{figure}{section}

\newtheorem{theorem}{Theorem}[section]

\newtheorem{lem}{Lemma}[section]
\newtheorem{prop}[theorem]{Proposition}
\newtheorem{defn}{Definition}[section]
\newtheorem{ex}{Example}[section]
\newtheorem{remark}{Remark}[section]

\makeatletter                 
\providecommand*{\cupdot}{%
	\mathbin{%
		\mathpalette\@cupdot{}%
	}%
}
\newcommand*{\@cupdot}[2]{%
	\ooalign{%
		$\m@th#1\cup$\cr
		\hidewidth$\m@th#1\cdot$\hidewidth
	}%
}
\makeatother

\makeatletter                  
\providecommand*{\bigcupdot}{%
	\mathbin{%
		\mathpalette\@bigcupdot{}%
	}%
}
\newcommand*{\@bigcupdot}[2]{%
	\ooalign{%
		$\m@th#1\bigcup$\cr
		\hidewidth$\m@th#1\cdot$\hidewidth
	}%
}
\makeatother

\title{Synchronous and asynchronous cyclic contractions on metric spaces}

\author{M\u ad\u alina P\u acurar$^{1}$}

\address{$^{1}$
	\indent Faculty of Economics and Business Administration\newline
	\indent Babe\c s-Bolyai University of Cluj-Napoca\newline 
	\indent T. Mihali 58-60, 400591 Cluj-Napoca, Romania}
\email{madalina.pacurar@econ.ubbcluj.ro}

\begin{document}
\maketitle \pagestyle{myheadings} \markboth{M. P\u ACURAR} {Synchronous and asynchronous cyclic contractions on metric spaces}

\subjclass[2020]{47H10; 54H25; 03Exx}

\keywords{\em cyclic covering; $r$-cyclic operator; synchronous $r$-cyclic contraction; asynchronous $r$-cyclic contraction; fixed point; weakly Picard operator; Picard operator; fixed point partition; iterative method; best proximity point}

\begin{abstract}
Motivated by the existence of cyclic phenomena in which some characteristics are mapped into corresponding ones over more than one phase, we introduce the $r$-cyclic operators with respect to a covering of a metric space and investigate their behavior. We study the convergence of the Picard iteration to a fixed point of such an operator under different types of generalized contraction conditions. The obtained results may have interesting practical applications in various research areas.
\end{abstract}

\section{Introduction} 

\medskip
In \cite{KirkSrinVeer} an interesting type of contractive condition for operators defined on metric spaces was suggested. The authors investigated there corresponding generalizations of important fixed point results, such as Banach's, Edelstein's or Caristi's fixed point theorems. They developed this way a technique which was used in the same paper in proving a known fixed point result for nonexpansive mappings.

This paper gave rise to a series of papers on what we are now used to call \textit{cyclic operators}. To our best knowledge, the notion of \textit{cyclic representation} first appeared in the short but interesting paper \cite{Rus-CyclicAnn}, which was also inspired by \cite{KirkSrinVeer}.

Since then many classes of cyclic operators have been studied, generally by extending known classes of generalized contractions on metric or Banach spaces (see for example \cite{Bilgili.Erhan.Karapinar.Turkoglu}, \cite{Harjani.Lopez.Sad-CyclicOrder}-\cite{Kim}, \cite{Kumari.Panthi}-\cite{Pacurar-CyclicBerinde}, 
\cite{PetruselG-CyclicPeriodic}, \cite{Weng.Liu.Chao}, for a very short list of them).

There are several other research areas where different types of cyclic phenomena are observed and studied. Of course, only some of these phenomena are likely to be approached by means of the instruments offered by our research. Here we mention only a few of these inspiring topics, where our results could be refined and adapted in order count as instruments in further studies and for developing new models: different types of oscillators (chemical, physical and others), social cycles, business cycles and waves, politico-economic cycles, financial cycles, ecosystems cycles, cycles in astrophysics and so on (see for example \cite{Christensen}, \cite{Foster.Shugart.Shuman}, \cite{Korotayev}, \cite{Wiescher}).

Inspired by these phenomena, we aim to introduce and study in the present paper cyclic contractive type operators corresponding to phenomena in which some characteristics do not appear in each and every phase of a cycle, but they appear with a certain periodicity (of $r$ phases). Besides new definitions and results that are given in the following, several remarks will complete the approach, mentioning aspects that can be relevant in applications. 

These new $r$-cyclic operators are an answer to natural problems and they prove to complete the family of cyclic type operators studied in the previous literature. They open this way a very wide space for further research, in order to extend the already existing results and theories referring to cyclic operators, if possible, in the case of $r$-cyclic operators. 

The \emph{Introduction} is followed by seven sections:

	\begin{itemize}
		\item[2.] \textit{Preliminaries} 
		\item[3.] \textit{Nuances of a definition. $r$-Cyclic operators}
		\item[4.] \textit{Properties of $r$-cyclic operators}
		\item[5.] \textit{Synchronous $r$-cyclic contractions on a metric space}
		\item[6.] \textit{Some remarks on synchronous $r$-cyclic contractions}
		\item[7.] \textit{Asynchronous $r$-cyclic contractions}
		\item[8.] \textit{Instead of a conclusion. The cyclic operators that have been missing}
	\end{itemize}

This first study of $r$-cyclic contractions does not exhaust the topic, as there are more questions to be answered and problems to be solved.

\section{Preliminaries} 
\medskip

A first (to our best knowledge) cyclic-type generalization of a classical fixed point result is Theorem 1.3 in \cite {KirkSrinVeer}, which we recall below in its original notation.

\begin{theorem}[\cite{KirkSrinVeer}]\label{Th_Kirk}
	Let $\{A_i\}_{i=1}^{p}$ be nonempty closed subsets of a complete metric space and suppose $F:\overset{p}{\underset{i=1}{\cup}}A_i\to\overset{p}{\underset{i=1}{\cup}}A_i$ satisfies the following conditions (where $A_{p+1}=A_1$):
	\begin{itemize}
		\item[(1)] $f(A_i)\subseteq A_{i+1}$, for $1\leq i \leq p$;
		\item[(2)] $\exists k\in(0,1)$ such that $d(F(x),F(y))\leq k\cdot d(x,y)$, $\forall x\in A_i, y\in A_{i+1}$, for $1\leq i \leq p$. 
	\end{itemize}
Then $F$ has a unique fixed point.
\end{theorem} 

The "secret" behind this result is that the Picard iteration associated to $F$ has infinitely many terms in each set $A_i$, $1\leq i \leq p$, and it is also a Cauchy sequence due to the contraction condition, these two arguments together with the completeness of the space and the closeness of $A_i,i=\overline{1,p}$ leading to the conclusion.

In \cite{Rus-CyclicAnn} the \textit{cyclic representation} relative to an operator $f:X\to X$ is defined as being $\overset{m}{\underset{i=1}{\bigcup}}X_i$ such that $X_i$ are all nonempty ($1\leq i \leq m$) and besides $f(X_1)\subseteq X_2, f(X_2)\subseteq X_3, \dots, f(x_{m-1})\subseteq X_m, f(X_m)\subseteq X_1$. Some examples of cyclic representations are given in \cite{Rus-CyclicAnn}. Other interesting examples of cyclic operators are to be found in the paper \cite{Horvat.Petric}.

The study of cyclic operators is still having an impressive development, judging after the number of papers in which the existence of fixed points or best proximity points is investigated, for various classes of operators and in various framework spaces.

This is why we found it interesting to revisit the papers that stood as starting point for this rich literature, to reconsider some notations and to reformulate some basic notions, having in view the phenomena they should be able to speak about, in mathematical language. This apparently trivial approach can lead to surprisingly interesting observations.

\medskip
\section{Nuances of a definition. $r$-Cyclic operators}
\medskip

In the rest of the paper we shall use the terminology dictated by the next definition, introduced in \cite{Rus-Set-Th.CJM}, which allows a more accurate language of cyclic operators and all related notions and results (instead of the initial \textit{cyclic representation} terminology).

\begin{defn}[\cite{Rus-Set-Th.CJM}]\label{Def_CycCov}
	Let $X$ be a nonempty set and $f:X \to X$ an operator. If there exists a covering of $X=\overset{m}{\underset{i=1}{\bigcup}}X_i$, $m\geq2$ such that
	\[
	f(X_1)\subseteq X_2, f(X_2)\subseteq X_3, \dots, f(x_{m-1})\subseteq X_m, f(X_m)\subseteq X_1,
	\] 
	then:
	\begin{itemize}
		\item[$i)$]  $\overset{m}{\underset{i=1}{\bigcup}}X_i$ is called a cyclic covering of $X$ w.r.t. $f$;
		\item[$ii)$] $f$ is called a cyclic operator w.r.t. the covering $\overset{m}{\underset{i=1}{\bigcup}}X_i$.
	\end{itemize}	
\end{defn}

\begin{remark} Note that in several papers (see for example \cite{Magadevan.Karpagam.Karapinar-pCyclic}) the term $p$-cyclic is used to indicate an operator that is cyclic in terms of Definition \ref{Def_CycCov}, where $p$ denotes the number of sets in the cyclic covering. This is still essentially different from what will be defined in the present paper as $r$-cyclic.
\end{remark}

\begin{remark}
	From the above definition one can see that $i)$ and $ii)$ are equivalent. Therefore in the sequel when one of them is mentioned the other one will be automatically assumed.
\end{remark}

The following simple example, included also in \cite{Pacurar.Rus-FPTCyclicOp}, will lead to an important remark regarding the notation used in the above definition.

\begin{ex}\label{Ex_X1X2X3}
	Let $X=X_1\cup X_2\cup X_3$ be a cyclic covering of $X$ w.r.t. the operator $f:X\to X$.
	
	One can easily check that according to the definition $X_2 \cup X_1\cup X_3$ is generally not a cyclic covering w.r.t. $f$, neither are $X_3\cup X_2\cup X_1$ or $X_1\cup X_3\cup X_2$.
	
	Still $X_2\cup X_3\cup X_1$ and $X_3\cup X_1\cup X_2$ are cyclic coverings w.r.t. $f$.
\end{ex}

\begin{remark}\label{Rem_CupCdot}
	It is now clear that the order in which the sets appear in the covering plays an important role, therefore none of them can change places. As the usual union of sets is commutative, we find that the notation $\overset{m}{\underset{i=1}{\bigcup}}X_i$ is not quite proper for a cyclic covering.
	
	So we introduce the notation 
	\[\overset{m}{\underset{i=1}{\cupdot}}X_i=X_1\cupdot X_2\cupdot\dots\cupdot X_m\]
	to indicate a cyclic covering w.r.t. to an operator, which will actually say that each of the $m$ cyclic permutations
	\begin{align*}
	& X_1\cup X_2 \cup X_3\cup\dots\cup X_{m-1}\cup X_m,\\
	& X_2\cup X_3 \cup\dots\cup X_{m-1}\cup X_m \cup X_1,\\
	&\vdots\\
	& X_m\cup X_1 \cup X_2 \cup\dots\cup X_{m-1} 
	\end{align*}
 is a cyclic covering of $X$ w.r.t. $f$, while generally any other permutation is not.      
\end{remark}

Now let us take a look at another simple but interesting example.

\begin{ex}
	If $\overset{m}{\underset{i=1}{\bigcupdot}}X_i$ is a cyclic covering w.r.t. $f:X\to X$ and $k\geq2$ is integer, then $f^k:X\to X$ is generally not a cyclic operator w.r.t. $\overset{m}{\underset{i=1}{\bigcupdot}}X_i$, but $f^{km+1}:X\to X$ is a cyclic operator w.r.t. $\overset{m}{\underset{i=1}{\bigcupdot}}X_i$.
\end{ex}

Motivated by the above observations and by the existence of cyclic phenomena where some studied characteristics or parameters are mapped not in the next coming phase, but in a more "remote" one (e.g., not from parents to children, but from grandparents to grandchildren, so over two generations), we introduce the following definition:

\begin{defn} \label{Def_rCycCov}
	Let $X$ be a nonempty set, $f:X \to X$ an operator and  $m\geq2$, $1\leq r\leq m$ integers. If there is a covering $X=\overset{m}{\underset{i=1}{\bigcupdot}}X_i$ such that
	\[
	f(X_1)\subseteq X_{1+r}, f(X_2)\subseteq X_{2+r}, \dots, f(X_m)\subseteq X_{m+r},
	\] 
	where for $p>m$ by $X_p$ we mean $X_{p\mod m}$, then:
\begin{itemize}
	\item[$i)$]  $\overset{m}{\underset{i=1}{\bigcupdot}}X_i$ is called a $r$-cyclic covering of $X$ w.r.t. $f$;
	\item[$ii)$] $f$ is called a $r$-cyclic operator w.r.t. the covering $\overset{m}{\underset{i=1}{\bigcupdot}}X_i$.
\end{itemize}
\end{defn}

\begin{remark}
	One can see that a cyclic covering or cyclic operator in the sense of Definition \ref{Def_CycCov} is a $1$-cyclic covering, respectively $1$-cyclic operator in the sense of Definition \ref{Def_rCycCov}. 
	
	For $r> m$, any $r$-cyclic operator/covering is actually $(r\mod m)$-cyclic operator/covering, that is why in the above definition we only consider $1\leq r \leq m$.
	
\end{remark}

\begin{remark}
	One can also check that the notation $\overset{m}{\underset{i=1}{\bigcupdot}}X_i$ from Remark \ref{Rem_CupCdot} remains consistent in the case of $r$-cyclic coverings as well, since if $\overset{m}{\underset{i=1}{\bigcupdot}}X_i$ is a $r$-cyclic covering w.r.t. an operator, then implicitly 
	\begin{align*}
	& X_1\cup X_2 \cup X_3\cup\dots\cup X_{m-1}\cup X_m,\\	
	& X_2\cup X_3 \cup\dots\cup X_{m-1}\cup X_m \cup X_1,\\	
	& \vdots\\
	& X_m\cup X_1 \cup X_2 \cup\dots\cup X_{m-1} 
	\end{align*}
	are $r$-cyclic w.r.t. to that operator, as well.
\end{remark}

Now let us analyze some simple examples, which are extensions of some examples in \cite{Rus-CyclicAnn}.

\begin{ex}
	If $\overset{m}{\underset{i=1}{\bigcupdot}}X_i$, $m\geq2$ is a cyclic covering w.r.t. $f:X\to X$, then it is a $k$-cyclic covering w.r.t. $f^k:X\to X$, for any $k\geq2$.
	
	See for example the case $r=2$. Since $f$ is cyclic w.r.t. $\overset{m}{\underset{i=1}{\bigcupdot}}X_i$, we have that 
	\[
	f^2(X_1)\subseteq f(f(X_1))\subseteq f(X_2)\subseteq X_3 \text{ etc.}	
	\]
	so 
	\[
	f^2(X_1)\subseteq X_3, f^2(X_2)\subseteq X_4, \dots, f^2(X_{m-2})\subseteq X_m,f^2(X_{m-1})\subseteq X_1,f^2(X_{m})\subseteq X_2,
	\]
	which by Definition \ref{Def_rCycCov} means exactly that $f^2$ is a $2$-cyclic operator w.r.t. $\overset{m}{\underset{i=1}{\bigcupdot}}X_i$.
\end{ex}

\begin{ex}
	Let $X=\{x_1, x_2, x_3,\dots, x_m\}$ a set, $X_i=\{x_i\},i=\overline{1,m}$ and $f:X\to X$ defined by 
	\begin{align*}
	&f(x_1)=x_4, f(x_2)=x_5, \dots, f(x_{m-3})=f(x_m),\\
	&f(x_{m-2})=f(x_1), f(x_{m-1})=f(x_2), f(x_{m})=f(x_3).
	\end{align*}
	Then $f$ is not a cyclic operator w.r.t. $\overset{m}{\underset{i=1}{\bigcupdot}}X_i$, but it is a $3$-cyclic operator w.r.t. to the same covering. 
	
	Note that it is generally not possible to "rearrange" the covering so that $f$ would be cyclic w.r.t. this rearrangement.
\end{ex}

\begin{ex}
	Let $X$ be a nonempty set and $f:X\to X$ an operator. Let $Y\subseteq X$ a nonempty strict subset of $X$ and let $n_0>0$, $1\leq r\leq m$ integers such that
	\[
	  X=Y\cup f^r(Y) \cup f^{2r}(Y) \cup \dots \cup f^{n_0r}(Y),
	\]
	with $f^{(n_0+1)r}(Y)\subseteq Y$.
	
	Then $f$ is a $r$-cyclic operator w.r.t. the covering $\overset{n_0}{\underset{i=0}{\bigcupdot}} f^{ir}(Y)$.
\end{ex}

\begin{remark}
	Generally if $X=\overset{m}{\underset{i=1}{\bigcupdot}}X_i$, $m\geq2$ is $r$-cyclic covering w.r.t. $f:X\to X$, then each $X_i$, $1\leq i\leq m$ is invariant for $f^m$.
	
	If in particular $m=k\cdot r$, with $k$ integer, then each $X_i$, $1\leq i\leq m$ is invariant for $f^k$.
\end{remark}

\begin{ex}\label{Ex_X1X2X3bis}
	We saw in Example \ref{Ex_X1X2X3} that if  $X=\overset{3}{\underset{i=1}{\bigcupdot}}X_i$ is a cyclic covering w.r.t. $f:X\to X$, then $X_1\cup X_2\cup X_3$, $X_2\cup X_3\cup X_1$ and $X_3\cup X_1\cup X_2$ are all cyclic coverings, while the remaining three cyclic permutations $X_1\cup X_3\cup X_2$, $X_2\cup X_1\cup X_3$ and $X_3\cup X_2\cup X_1$ are generally not cyclic coverings.
	
	Still one can check that they are all $2$-cyclic coverings w.r.t. $f$.
\end{ex}

\bigskip
\section{Properties of $r$-cyclic operators}
\bigskip

The above Example \ref{Ex_X1X2X3bis} could lead us to the false conclusion that, given an operator $f:\overset{m}{\underset{i=1}{\bigcupdot}}X_i\to \overset{m}{\underset{i=1}{\bigcupdot}}X_i$, no matter how we rearrange the elements of the covering $\overset{m}{\underset{i=1}{\bigcupdot}}X_i$, it is possible to find a positive integer $r$ such that $f$ is $r$-cyclic w.r.t. to this rearrangement. 

This is generally not true, and an answer will be given by Lemmas \ref{Lem_mrk} and \ref{Lem_mr1}. These results can be simply deduced by observing the behavior of $r$-cyclic coverings w.r.t. $f$, for various values of $m$ and $r$.

Before stating the lemmas, we propose a visualization of some $r$-cyclic coverings, as we found this approach extremely useful for a better understanding of the results to come. Note that the sets $X_1, X_2, \dots, X_m$ are not supposed to be disjoint, in general.

According to Definition \ref{Def_rCycCov}, we refer to coverings $X=\overset{m}{\underset{i=1}{\bigcupdot}}X_i$, $m\geq2$ and to an operator $f:X\to X$ such that
	\[
	f(X_i)\subseteq X_{i+r}, i=\overline{1,m},
	\]
for $1\leq r\leq m$, where $X_p=X_{p\mod m}$ for any $p>m$.

\subsection{The case $m=2$}\ 

Then $X=X_1\cupdot X_2$ and we have to analyze two values of $r$.

For $r=1$, $f$ would be $1$-cyclic (or simply cyclic), that is, \[f(X_1)\subseteq X_2, f(X_2)\subseteq X_1.\]

Graphically we shall represent this as

\begin{figure}[!h]
	\begin{tikzpicture}[>=stealth]
		\node (1) at (180:1) {$X_1$};
		\node (2) at (0:1) {$X_2$};
		\draw [myblue,->] (1) to  [bend right, looseness=1] (2);
		\draw [myblue,->] (2) to  [bend right, looseness=1] (1);
	\end{tikzpicture}
	\caption{Cyclic covering for $m=2,r=1$}
\end{figure}
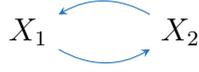

For $r=2$, $f$ would be $2$-cyclic, that is, \[f(X_1)\subseteq X_1, f(X_2)\subseteq X_2.\] We shall represent this as

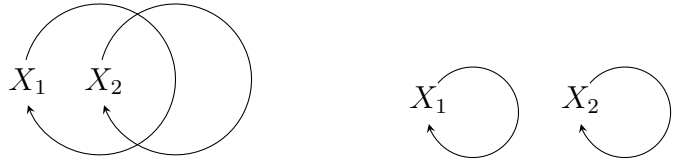
\begin{figure}[!h]
	\begin{tikzpicture}[>=stealth]
	
		\draw[black,->] (0,0) arc [start angle = 165, end angle =-160, radius = 10mm] node [anchor=south] {\color{black}{$X_1$}}; 
		\draw[black,->] (1,0) arc [start angle = 165, end angle =-160, radius = 10mm] node [anchor=south] {\color{black}{$X_2$}};	 
	\end{tikzpicture}
\hspace{0.5cm}
\ \ \ 
\hspace{0.5cm}
	\begin{tikzpicture}[>=stealth]
		\draw[black,->] (1) arc [start angle = 140, end angle =-165, radius = 6mm] node [anchor=south] {\color{black}{$X_1$}}; 
		\draw[black,->] (2) arc [start angle = 140, end angle =-165, radius = 6mm] node [anchor=south] {\color{black}{$X_2$}};	 
	\end{tikzpicture}
	\caption{Cyclic covering for $m=2,r=2$, direct vs. simplified}
\end{figure}

\subsection{The case $m=3$}\ 

Then $X=X_1\cupdot X_2\cupdot X_3$ and we have to analyze three values of $r$.

For $r=1$, $f$ would be $1$-cyclic (or simply cyclic), that is, \[f(X_1)\subseteq X_2, f(X_2)\subseteq X_3, f(X_3)\subseteq X_1.\]

Graphically we shall represent this as

\begin{figure}[!h]
	\begin{tikzpicture}[>=stealth]
		\foreach \x in {1,...,3}
		\node (\x) at (180-\x*360/3:1) {$X_{\x}$};	
		\draw [myblue,->] (1) to  [bend left, looseness=1] (2);
		\draw [myblue,->] (2) to  [bend left, looseness=1] (3);
		\draw [myblue,->] (3) to  [bend left, looseness=1] (1);	
	\end{tikzpicture}
	\caption{Cyclic covering for $m=3,r=1$}
\end{figure}
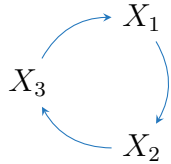

For $r=2$, $f$ would be $2$-cyclic, that is, \[f(X_1)\subseteq X_3, f(X_2)\subseteq X_1, f(X_3)\subseteq X_2.\]

Graphically we shall represent this as

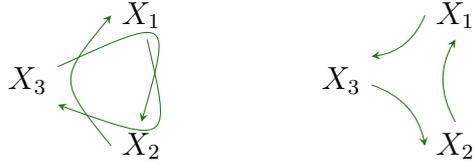
\begin{figure}[!h]
	\begin{tikzpicture}[>=stealth]
		\foreach \x in {1,...,3}
		\node (\x) at (180-\x*360/3:1) {$X_{\x}$};	
		\draw [mygreen,->] (1) .. controls (2.1) .. (3.south east);
		\draw [mygreen,->] (2.west) .. controls (3.1) .. (1.west);
		\draw [mygreen,->] (3) .. controls (1.1) .. (2.north);	
	\end{tikzpicture}
\hspace{0.5cm}
\ \ \ 
\hspace{0.5cm}
	\begin{tikzpicture}[>=stealth] 
		\foreach \x in {1,...,3}
		\node (\x) at (180-\x*360/3:1) {$X_{\x}$};	
		\draw [mygreen,->] (1.west) to  [bend left, looseness=1] (3.north east);
		\draw [mygreen,->] (2.north) to  [bend left, looseness=1] (1.south);
		\draw [mygreen,->] (3) to  [bend left, looseness=1] (2.west);	
	\end{tikzpicture}
\caption{Cyclic covering for $m=3,r=2$, direct vs. simplified}	
\end{figure}

Note that in the following we shall adopt this kind of "left bent" arrows like in the second representation to tell that $f$ maps $X_1$ into $X_3$, stepping "over" $X_2$, $X_2$ into $X_1$ and so on.

For $r=3$ we obviously have

\begin{figure}[H]
	\begin{tikzpicture}[>=stealth]
		\draw[->] (1) arc [start angle = 150, end angle =-165, radius = 7mm] node [anchor=south] {\color{black}{$X_1$}}; 
		\draw[->] (2) arc [start angle = 150, end angle =-165, radius = 7mm] node [anchor=south] {\color{black}{$X_2$}};
		\draw[->] (3) arc [start angle = 150, end angle =-165, radius = 7mm] node [anchor=south] {\color{black}{$X_3$}};
	\end{tikzpicture}
\caption{Cyclic covering for $m=3,r=3$}	
\end{figure}
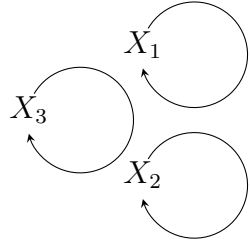

\subsection{The case $m=4$}\ 

Now $X=X_1\cupdot X_2\cupdot X_3\cupdot X_4$ and we have to analyze four values of $r$, depicted in the next  Figures \ref{Fig m4-r1.2} and \ref{Fig m4-r3.4} below.

\begin{figure}[!h]
	\begin{tikzpicture}[>=stealth,scale=1.4] 
		\foreach \x in {1,...,4}
		\node (\x) at (180-\x*360/4:1) {$X_{\x}$};	
		\draw [myblue,->] (1) to  [bend left, looseness=1] (2);
		\draw [myblue,->] (2) to  [bend left, looseness=1] (3);
		\draw [myblue,->] (3) to  [bend left, looseness=1] (4);	
		\draw [myblue,->] (4) to  [bend left, looseness=1] (1);	
	\end{tikzpicture}
	\hspace{1.5cm}
	\begin{tikzpicture}[>=stealth,scale=1.5] 
		\foreach \x in {1,...,4}
		\node (\x) at (180-\x*360/4:1) {$X_{\x}$};	
		\draw [plum,->] (1) to  [bend left, looseness=1] (3);
		\draw [tangelo,->] (2) to  [bend left, looseness=1] (4);
		\draw [plum,->] (3) to  [bend left, looseness=1] (1);	
		\draw [tangelo,->] (4) to  [bend left, looseness=1] (2);	
	\end{tikzpicture}
	\caption{Cyclic coverings for $m=4,r=1$ and $m=4,r=2$}
	\label{Fig m4-r1.2}	
\end{figure}
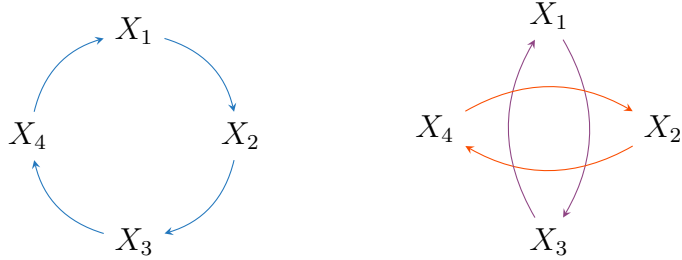

\begin{figure}[!h]
	\begin{tikzpicture}[>=stealth,scale=1.5] 
		\foreach \x in {1,...,4}
		\node (\x) at (180-\x*360/4:1) {$X_{\x}$};	
		\draw [mygreen,->] (1) to  [bend left, looseness=1] (4);
		\draw [mygreen,->] (2) to  [bend left, looseness=1] (1);
		\draw [mygreen,->] (3) to  [bend left, looseness=1] (2);	
		\draw [mygreen,->] (4) to  [bend left, looseness=1] (3);	
	\end{tikzpicture}
	\hspace{1.5cm}
	\begin{tikzpicture}[>=stealth]
		\draw[->] (1) arc [start angle = 145, end angle =-165, radius = 8mm] node [anchor=south] {\color{black}{$X_1$}}; 
		\draw[->] (2) arc [start angle = 145, end angle =-165, radius = 8mm] node [anchor=south] {\color{black}{$X_2$}};
		\draw[->] (3) arc [start angle = 145, end angle =-165, radius = 8mm] node [anchor=south] {\color{black}{$X_3$}};
		\draw[->] (4) arc [start angle = 145, end angle =-165, radius = 8mm] node [anchor=south] {\color{black}{$X_4$}};
	\end{tikzpicture}
	\caption{Cyclic coverings for $m=4,r=3$ and $m=4,r=4$}
	\label{Fig m4-r3.4}
\end{figure}
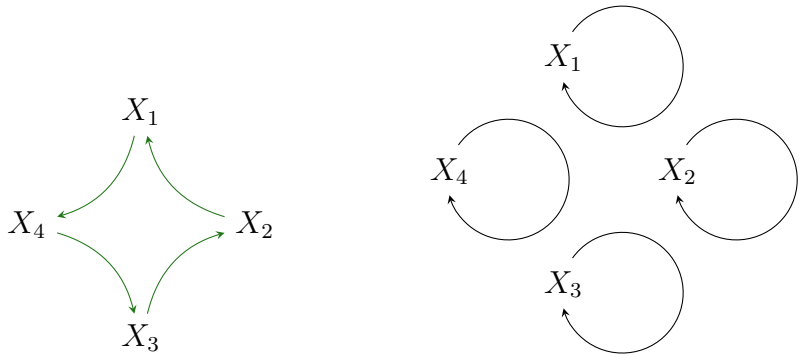

\subsection{Comments on visual representations}\ 

Although it might seem superfluous, it makes sense going further with this visualization approach, as it says more than a proof, especially because rigorous notation in the case of $r$-cyclic operators is getting sometimes really hairy, as we shall see later.

One can already notice that in some cases, no matter which of $X_i$, $i=\overline{1,m}$ is the starting set, the operator $f$ will run through each of the sets of the covering. In other cases, some "closed circuits" can be observed (for example when $m=4, r=2$), that is, if the starting point is in one of the sets contained in this circuit, then it will never be mapped outside the circuit.

\begin{remark}\label{Rem_circuit}
	In the rest of the paper we shall use this term of \textit{closed circuit} or simply \textit{circuit} in order to denote such a family of sets, invariant for $f$.
\end{remark}

We try to make this phenomenon more visible in our representations by means of the color code we use: blue for a $1$-cyclic operator, black for $m$-cyclic operators, green in all those cases when $f$ runs through all the sets $X_i$, $i=\overline{1,m}$. Finally, if there are closed circuits to be noticed, we use other colors, one different color for each such circuit. This color code can already be noticed for the case $m=4$ above. We still have to mention that all these representations can be well understood without any color, only by attentively following the arrows.

\subsection{The case $m=5$}\ 

Here $X=X_1\cupdot X_2\cupdot X_3\cupdot X_4\cupdot X_5$ and we have to analyze five different values of $r$.

For $r=1$, $r=2$, respectively $r=3$, the representations are
\begin{figure}[H]
	\begin{tikzpicture}[>=stealth,scale=1.5]  
		\foreach \x in {1,...,5}
		\node (\x) at (180-\x*360/5:1) {$X_{\x}$};	
		\draw [myblue,->] (1) to  [bend left, looseness=1] (2);
		\draw [myblue,->] (2) to  [bend left, looseness=1] (3);
		\draw [myblue,->] (3) to  [bend left, looseness=1] (4);	
		\draw [myblue,->] (4) to  [bend left, looseness=1] (5);	
		\draw [myblue,->] (5) to  [bend left, looseness=1] (1);		
	\end{tikzpicture}
	\hspace{0.25cm}
	\begin{tikzpicture}[>=stealth,scale=1.5]  
		\foreach \x in {1,...,5}
		\node (\x) at (180-\x*360/5:1) {$X_{\x}$};	
		\draw [mygreen,->] (1) to  [bend left, looseness=1] (3);
		\draw [mygreen,->] (2) to  [bend left, looseness=1] (4);
		\draw [mygreen,->] (3) to  [bend left, looseness=1] (5);	
		\draw [mygreen,->] (4) to  [bend left, looseness=1] (1);	
		\draw [mygreen,->] (5) to  [bend left, looseness=1] (2);
	\end{tikzpicture}
	\hspace{0.25cm}
	\begin{tikzpicture}[>=stealth,scale=1.5]  
		\foreach \x in {1,...,5}
		\node (\x) at (180-\x*360/5:1) {$X_{\x}$};	
		\draw [mygreen,->] (1) to  [bend left, looseness=1] (4);
		\draw [mygreen,->] (2) to  [bend left, looseness=1] (5);
		\draw [mygreen,->] (3) to  [bend left, looseness=1] (1);	
		\draw [mygreen,->] (4) to  [bend left, looseness=1] (2);	
		\draw [mygreen,->] (5) to  [bend left, looseness=1] (3);
	\end{tikzpicture}
	\caption{Cyclic coverings for $m=5$ and $r=1,r=2,r=3$}
\end{figure}
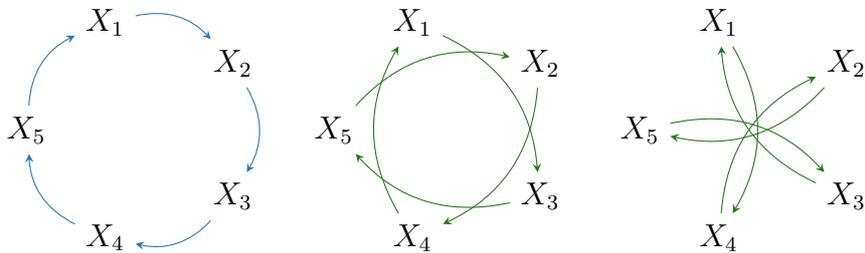

whereas for $r=4$ and, respectively, $r=5$ we have

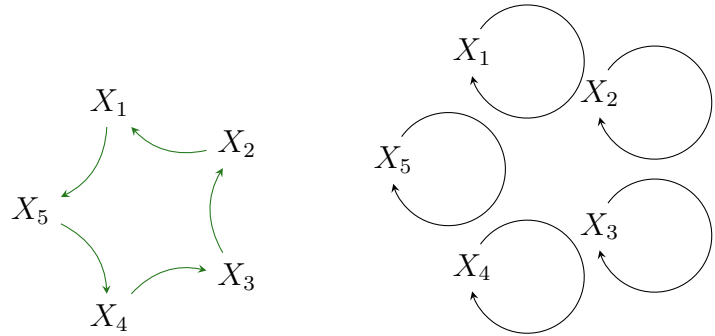
\begin{figure}[H]
	\begin{tikzpicture}[>=stealth,scale=1.5]  
		\foreach \x in {1,...,5}
		\node (\x) at (180-\x*360/5:1) {$X_{\x}$};	
		\draw [mygreen,->] (1) to  [bend left, looseness=1] (5);
		\draw [mygreen,->] (2) to  [bend left, looseness=1] (1);
		\draw [mygreen,->] (3) to  [bend left, looseness=1] (2);	
		\draw [mygreen,->] (4) to  [bend left, looseness=1] (3);	
		\draw [mygreen,->] (5) to  [bend left, looseness=1] (4);
	\end{tikzpicture}
	\hspace{1cm}
	\begin{tikzpicture}[>=stealth,scale=1.5]
		\draw[->] (1) arc [start angle = 145, end angle =-165, radius = 5mm] node [anchor=south] {\color{black}{$X_1$}}; 
		\draw[->] (2) arc [start angle = 145, end angle =-165, radius = 5mm] node [anchor=south] {\color{black}{$X_2$}};
		\draw[->] (3) arc [start angle = 145, end angle =-165, radius = 5mm] node [anchor=south] {\color{black}{$X_3$}};
		\draw[->] (4) arc [start angle = 145, end angle =-165, radius = 5mm] node [anchor=south] {\color{black}{$X_4$}};
		\draw[->] (5) arc [start angle = 145, end angle =-165, radius = 5mm] node [anchor=south] {\color{black}{$X_5$}};
	\end{tikzpicture}
	\caption{Cyclic coverings for $m=5$ and $r=4,r=5$}
\end{figure}

For all the cases $r=1,2,3,4$ the operator $f$ will run through all the sets of the covering, no matter which the starting set is. 

\subsection{The case $m=6$}\ 

Now $X=X_1\cupdot X_2\cupdot X_3\cupdot X_4\cupdot X_5\cupdot X_6$ and we have to analyze six different values of $r$. This case will show interesting behavior.

For $r=1$, $r=2$, respectively $r=3$, the representations are:
\begin{figure}[H]
	\begin{tikzpicture}[>=stealth,scale=1.5] 
		\foreach \x in {1,...,6}
		\node (\x) at (180-\x*360/6:1) {$X_{\x}$};	
		\draw [myblue,->] (1) to  [bend left, looseness=1] (2);
		\draw [myblue,->] (2) to  [bend left, looseness=1] (3);
		\draw [myblue,->] (3) to  [bend left, looseness=1] (4);	
		\draw [myblue,->] (4) to  [bend left, looseness=1] (5);	
		\draw [myblue,->] (5) to  [bend left, looseness=1] (6);	
		\draw [myblue,->] (6) to  [bend left, looseness=1] (1);		
	\end{tikzpicture}
	\hspace{0.25cm}
	\begin{tikzpicture}[>=stealth,scale=1.5] 
		\foreach \x in {1,...,6}
		\node (\x) at (180-\x*360/6:1) {$X_{\x}$};	
		\draw [plum,->] (1.south east) to  [bend left, looseness=1] (3.north west);
		\draw [tangelo,->] (2) to  [bend left, looseness=1] (4);
		\draw [plum,->] (3.west) to  [bend left, looseness=1] (5.north east);	
		\draw [tangelo,->] (4) to  [bend left, looseness=1] (6);	
		\draw [plum,->] (5.north) to  [bend left, looseness=1] (1.south);
		\draw [tangelo,->] (6) to  [bend left, looseness=0.8] (2);
	\end{tikzpicture}
	\hspace{0.25cm}
	\begin{tikzpicture}[>=stealth,scale=1.5] 
		\foreach \x in {1,...,6}
		\node (\x) at (180-\x*360/6:1) {$X_{\x}$};	
		\draw [plum,->] (1) to  [bend left, looseness=1] (4);
		\draw [tangelo,->] (2) to  [bend left, looseness=1] (5);
		\draw [yellow,->] (3) to  [bend left, looseness=1] (6);	
		\draw [plum,->] (4) to  [bend left, looseness=1] (1);	
		\draw [tangelo,->] (5) to  [bend left, looseness=1] (2);
		\draw [yellow,->] (6) to  [bend left, looseness=1] (3);
	\end{tikzpicture}
	\caption{Cyclic coverings for $m=6$ and $r=1,r=2,r=3$}
\end{figure}

whereas for $r=4$, $r=5$ and, respectively, $r=6$ we have:

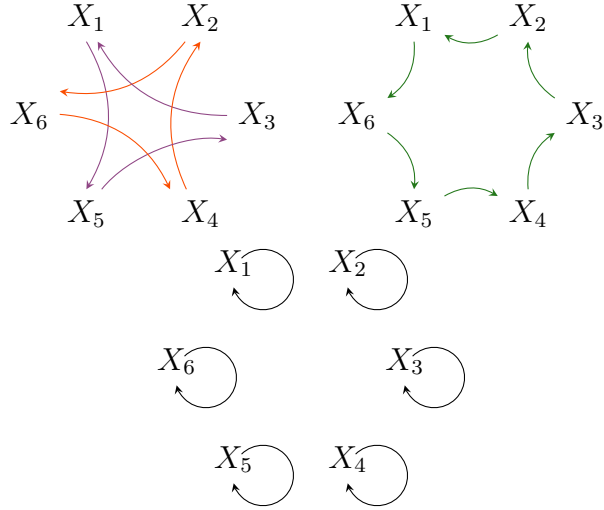
\begin{figure}[H]
	\begin{tikzpicture}[>=stealth,scale=1.5] 
		\foreach \x in {1,...,6}
		\node (\x) at (180-\x*360/6:1) {$X_{\x}$};	
		\draw [plum,->] (1.south) to  [bend left, looseness=1] (5.north);
		\draw [tangelo,->] (2) to  [bend left, looseness=1] (6.north east);
		\draw [plum,->] (3.west) to  [bend left, looseness=1] (1);	
		\draw [tangelo,->] (4) to  [bend left, looseness=1] (2.south);	
		\draw [plum,->] (5) to  [bend left, looseness=0.8] (3.south west);
		\draw [tangelo,->] (6) to  [bend left, looseness=1] (4.north west);
	\end{tikzpicture}
	\hspace{0.25cm}
	\begin{tikzpicture}[>=stealth,scale=1.5] 
		\foreach \x in {1,...,6}
		\node (\x) at (180-\x*360/6:1) {$X_{\x}$};	
		\draw [mygreen,->] (1) to  [bend left, looseness=1] (6);
		\draw [mygreen,->] (2) to  [bend left, looseness=1] (1);
		\draw [mygreen,->] (3) to  [bend left, looseness=1] (2);	
		\draw [mygreen,->] (4) to  [bend left, looseness=1] (3);	
		\draw [mygreen,->] (5) to  [bend left, looseness=1] (4);
		\draw [mygreen,->] (6) to  [bend left, looseness=1] (5);	
	\end{tikzpicture}
	\hspace{0.25cm}
	\begin{tikzpicture}[>=stealth]
		\draw[->] (1) arc [start angle = 135, end angle =-165, radius = 4mm] node [anchor=south] {\color{black}{$X_1$}}; 
		\draw[->] (2) arc [start angle = 135, end angle =-165, radius = 4mm] node [anchor=south] {\color{black}{$X_2$}};
		\draw[->] (3) arc [start angle = 135, end angle =-165, radius = 4mm] node [anchor=south] {\color{black}{$X_3$}};
		\draw[->] (4) arc [start angle = 135, end angle =-165, radius = 4mm] node [anchor=south] {\color{black}{$X_4$}};
		\draw[->] (5) arc [start angle = 135, end angle =-165, radius = 4mm] node [anchor=south] {\color{black}{$X_5$}};
		\draw[->] (6) arc [start angle = 135, end angle =-165, radius = 4mm] node [anchor=south] {\color{black}{$X_6$}};
	\end{tikzpicture}
	\caption{Cyclic coverings for $m=6$ and $r=4,r=5,r=6$}
\end{figure}

As one can easily notice, the only cases when $f$ runs through all the sets of the covering are when $m=1$ and $m=5$. For $m=2$ and $m=4$ there appear two closed circuits and for $m=3$ there are three such circuits.

Note that, although the involved sets are the same, the circuits do not coincide for $m=2$ ($f$ is 2-cyclic) and $m=4$ ($f$ is 4-cyclic), since 
	\[ X_1\cupdot X_3\cupdot X_5 \text{ is not the same as } X_1\cupdot X_5\cupdot X_3	\]
 and 
 	\[ X_2\cupdot X_4\cupdot X_6 \text{ is not the same as } X_2\cupdot X_6\cupdot X_4.	\]

\subsection{The case $m=10$}\ 

	Though the other cases are not less interesting, we omit them and we represent below only the resulting configurations when $m=10$ and $r=1,2,\dots,10$. In this case the cyclic covering is $X=X_1\cupdot X_2\cupdot \dots \cupdot X_{10}$.

	For $r=1$ and, respectively, $r=2$, the representations are:
	\begin{figure}[H]
		\begin{tikzpicture}[>=stealth,scale=2]
			\foreach \x in {1,...,10}
			\node (\x) at (90-\x*360/10:1) {$X_{\x}$};
			\draw [myblue,->] (1) to  [bend left, looseness=1] (2);
			\draw [myblue,->] (2) to  [bend left, looseness=1] (3);
			\draw [myblue,->] (3) to  [bend left, looseness=1] (4);	
			\draw [myblue,->] (4) to  [bend left, looseness=1] (5);	
			\draw [myblue,->] (5) to  [bend left, looseness=1] (6);	
			\draw [myblue,->] (6) to  [bend left, looseness=1] (7);	
			\draw [myblue,->] (7) to  [bend left, looseness=1] (8);
			\draw [myblue,->] (8) to  [bend left, looseness=1] (9);	
			\draw [myblue,->] (9) to  [bend left, looseness=1] (10);	
			\draw [myblue,->] (10) to  [bend left, looseness=1] (1);	
		\end{tikzpicture}
		\hspace{1cm}
		\begin{tikzpicture}[>=stealth,scale=2]    
			\foreach \x in {1,...,10}
			\node (\x) at (90-\x*360/10:1) {$X_{\x}$};
			\draw [plum,->] (1.east) to  [bend left, looseness=1] (3.east);
			\draw [tangelo,->] (2) to  [bend left, looseness=0.2] (4.north);
			\draw [plum,->] (3.south) to  [bend left, looseness=1] (5.south east);	
			\draw [tangelo,->] (4.west) to  [bend left, looseness=0.2] (6.east);	
			\draw [plum,->] (5.south west) to  [bend left, looseness=1] (7.south);	
			\draw [tangelo,->] (6.north) to  [bend left, looseness=0.2] (8.south);	
			\draw [plum,->] (7.west) to  [bend left, looseness=1] (9.west);
			\draw [tangelo,->] (8.north east) to  [bend left, looseness=0.2] (10.south west);	
			\draw [plum,->] (9.north) to  [bend left, looseness=1] (1.north);	
			\draw [tangelo,->] (10.south east) to  [bend left, looseness=0.2] (2.north west);	
		\end{tikzpicture}
		\caption{Cyclic coverings for $m=10$ and $r=1,r=2$}
	\end{figure}
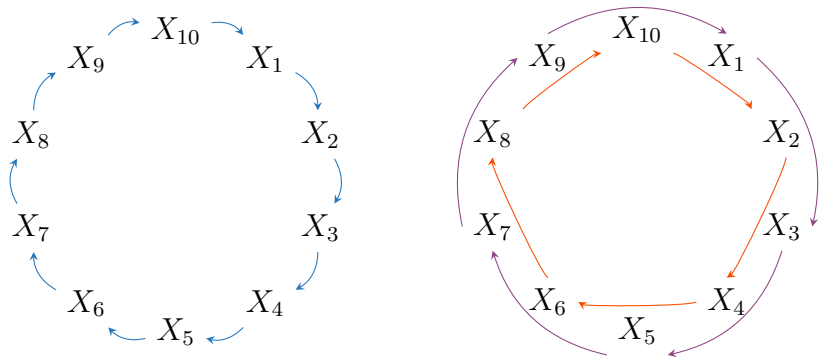

	For $r=3$ and, respectively, $r=4$, one has:
	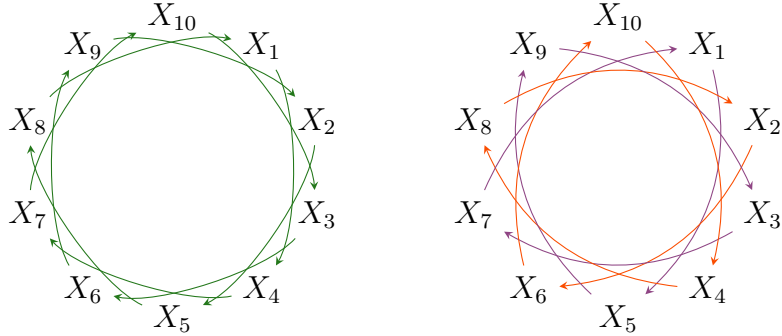
\begin{figure}[H]
		\begin{tikzpicture}[>=stealth,scale=2] 
		\foreach \x in {1,...,10}
		\node (\x) at (90-\x*360/10:1) {$X_{\x}$};
		\draw [mygreen,->] (1) to  [bend left, looseness=0.6] (4);
		\draw [mygreen,->] (2) to  [bend left, looseness=0.6] (5);
		\draw [mygreen,->] (3) to  [bend left, looseness=0.6] (6);	
		\draw [mygreen,->] (4) to  [bend left, looseness=0.6] (7);	
		\draw [mygreen,->] (5) to  [bend left, looseness=0.6] (8);	
		\draw [mygreen,->] (6) to  [bend left, looseness=0.6] (9);	
		\draw [mygreen,->] (7) to  [bend left, looseness=0.6] (10);
		\draw [mygreen,->] (8) to  [bend left, looseness=0.6] (1);	
		\draw [mygreen,->] (9) to  [bend left, looseness=0.6] (2);	
		\draw [mygreen,->] (10) to  [bend left, looseness=0.6] (3);	
		\end{tikzpicture}	
	\hspace{1cm}
		\begin{tikzpicture}[>=stealth,scale=2]  
			\foreach \x in {1,...,10}
			\node (\x) at (90-\x*360/10:1) {$X_{\x}$};
			\draw [plum,->] (1) to  [bend left, looseness=1] (5);
			\draw [tangelo,->] (2) to  [bend left, looseness=1] (6);
			\draw [plum,->] (3) to  [bend left, looseness=1] (7);	
			\draw [tangelo,->] (4) to  [bend left, looseness=1] (8);	
			\draw [plum,->] (5) to  [bend left, looseness=1] (9);	
			\draw [tangelo,->] (6) to  [bend left, looseness=1] (10);	
			\draw [plum,->] (7) to  [bend left, looseness=1] (1);
			\draw [tangelo,->] (8) to  [bend left, looseness=1] (2);	
			\draw [plum,->] (9) to  [bend left, looseness=1] (3);	
			\draw [tangelo,->] (10) to  [bend left, looseness=1] (4);	
		\end{tikzpicture}
		\caption{Cyclic coverings for $m=10$ and $r=3,r=4$}
	\end{figure}

	For $r=5$ and, respectively, $r=6$, the representations are:
	\begin{figure}[H]
		\begin{tikzpicture}[>=stealth,scale=2]
			\foreach \x in {1,...,10}
			\node (\x) at (90-\x*360/10:1) {$X_{\x}$};
			\draw [myblue,->] (1) to  [bend left, looseness=1] (6);
			\draw [yellow,->] (2) to  [bend left, looseness=1] (7);
			\draw [mygreen,->] (3) to  [bend left, looseness=1] (8);	
			\draw [tangelo,->] (4) to  [bend left, looseness=1] (9);	
			\draw [plum,->] (5) to  [bend left, looseness=1] (10);	
			\draw [myblue,->] (6) to  [bend left, looseness=1] (1);	
			\draw [yellow,->] (7) to  [bend left, looseness=1] (2);
			\draw [mygreen,->] (8) to  [bend left, looseness=1] (3);	
			\draw [tangelo,->] (9) to  [bend left, looseness=1] (4);	
			\draw [plum,->] (10) to  [bend left, looseness=1] (5);	
		\end{tikzpicture}	
		\hspace{1cm}
		\begin{tikzpicture}[>=stealth,scale=2] 
			\foreach \x in {1,...,10}
			\node (\x) at (90-\x*360/10:1) {$X_{\x}$};
			\draw [plum,->] (1) to  [bend left, looseness=0.5] (7);
			\draw [tangelo,->] (2) to  [bend left, looseness=1] (8);
			\draw [plum,->] (3) to  [bend left, looseness=0.5] (9);	
			\draw [tangelo,->] (4) to  [bend left, looseness=1] (10);	
			\draw [plum,->] (5) to  [bend left, looseness=0.5] (1);	
			\draw [tangelo,->] (6) to  [bend left, looseness=1] (2);	
			\draw [plum,->] (7) to  [bend left, looseness=0.5] (3);
			\draw [tangelo,->] (8) to  [bend left, looseness=1] (4);	
			\draw [plum,->] (9) to  [bend left, looseness=0.5] (5);	
			\draw [tangelo,->] (10) to  [bend left, looseness=1] (6);	
		\end{tikzpicture}
		\caption{Cyclic coverings for $m=10$ and $r=5,r=6$}
	\end{figure}
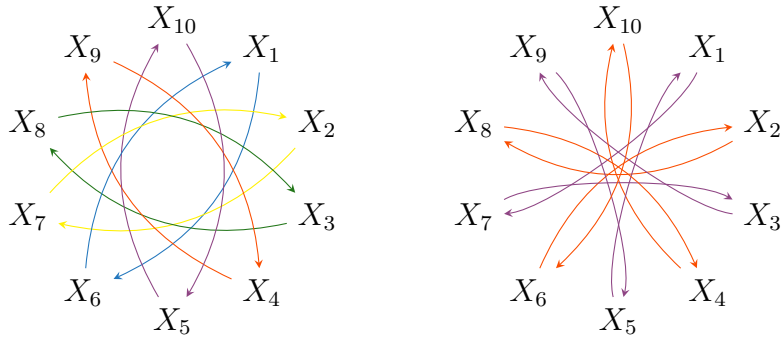

	For $r=7$, respectively $r=8$, one obtains:
	
	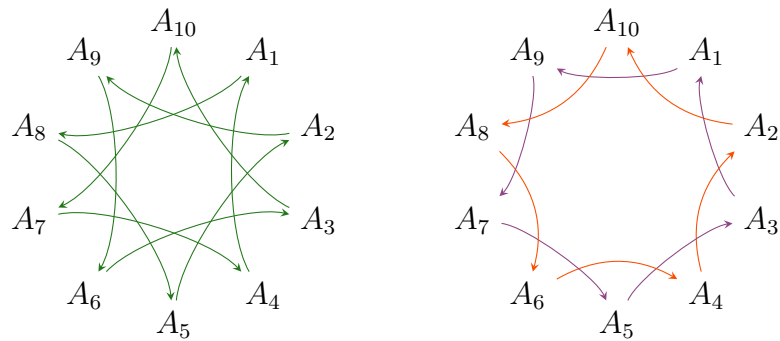
\begin{figure}[H]
		\begin{tikzpicture}[>=stealth,scale=2] 
			\foreach \x in {1,...,10}
			\node (\x) at (90-\x*360/10:1) {$A_{\x}$};
			\draw [mygreen,->] (1) to  [bend left, looseness=0.6] (8);
			\draw [mygreen,->] (2) to  [bend left, looseness=0.6] (9);
			\draw [mygreen,->] (3) to  [bend left, looseness=0.6] (10);	
			\draw [mygreen,->] (4) to  [bend left, looseness=0.6] (1);	
			\draw [mygreen,->] (5) to  [bend left, looseness=0.6] (2);	
			\draw [mygreen,->] (6) to  [bend left, looseness=0.6] (3);	
			\draw [mygreen,->] (7) to  [bend left, looseness=0.6] (4);
			\draw [mygreen,->] (8) to  [bend left, looseness=0.6] (5);	
			\draw [mygreen,->] (9) to  [bend left, looseness=0.6] (6);	
			\draw [mygreen,->] (10) to  [bend left, looseness=0.6] (7);	
		\end{tikzpicture}
		\hspace{1cm}
		\begin{tikzpicture}[>=stealth,scale=2] 
			\foreach \x in {1,...,10}
			\node (\x) at (90-\x*360/10:1) {$A_{\x}$};
			\draw [plum,->] (1) to  [bend left, looseness=0.5] (9);
			\draw [tangelo,->] (2) to  [bend left, looseness=1] (10);
			\draw [plum,->] (3) to  [bend left, looseness=0.5] (1);	
			\draw [tangelo,->] (4) to  [bend left, looseness=1] (2);	
			\draw [plum,->] (5) to  [bend left, looseness=0.5] (3);	
			\draw [tangelo,->] (6) to  [bend left, looseness=1] (4);	
			\draw [plum,->] (7) to  [bend left, looseness=0.5] (5);
			\draw [tangelo,->] (8) to  [bend left, looseness=1] (6);	
			\draw [plum,->] (9) to  [bend left, looseness=0.5] (7);	
			\draw [tangelo,->] (10) to  [bend left, looseness=1] (8);	
		\end{tikzpicture}
		\caption{Cyclic coverings for $m=10$ and $r=7,r=8$}
	\end{figure}

	And finally for $r=9$, respectively $r=10$, the representations are:
	
	\begin{figure}[H]
		\begin{tikzpicture}[>=stealth,scale=2] 
			\foreach \x in {1,...,10}
			\node (\x) at (90-\x*360/10:1) {$A_{\x}$};
			\draw [mygreen,->] (1) to  [bend left, looseness=1.5] (10);
			\draw [mygreen,->] (2) to  [bend left, looseness=1.5] (1);
			\draw [mygreen,->] (3) to  [bend left, looseness=1.5] (2);	
			\draw [mygreen,->] (4) to  [bend left, looseness=1.5] (3);	
			\draw [mygreen,->] (5) to  [bend left, looseness=1.5] (4);	
			\draw [mygreen,->] (6) to  [bend left, looseness=1.5] (5);	
			\draw [mygreen,->] (7) to  [bend left, looseness=1.5] (6);
			\draw [mygreen,->] (8) to  [bend left, looseness=1.5] (7);	
			\draw [mygreen,->] (9) to  [bend left, looseness=1.5] (8);	
			\draw [mygreen,->] (10) to  [bend left, looseness=1.5] (9);	
		\end{tikzpicture}
		\hspace{0.25cm}
	\begin{tikzpicture}[>=stealth]
		\draw[->] (1) arc [start angle = 135, end angle =-165, radius = 4mm] node [anchor=south] {\color{black}{$X_1$}}; 
		\draw[->] (2) arc [start angle = 135, end angle =-165, radius = 4mm] node [anchor=south] {\color{black}{$X_2$}};
		\draw[->] (3) arc [start angle = 135, end angle =-165, radius = 4mm] node [anchor=south] {\color{black}{$X_3$}};
		\draw[->] (4) arc [start angle = 135, end angle =-165, radius = 4mm] node [anchor=south] {\color{black}{$X_4$}};
		\draw[->] (5) arc [start angle = 135, end angle =-165, radius = 4mm] node [anchor=south] {\color{black}{$X_5$}};
		\draw[->] (6) arc [start angle = 135, end angle =-165, radius = 4mm] node [anchor=south] {\color{black}{$X_6$}};
		\draw[->] (7) arc [start angle = 135, end angle =-165, radius = 4mm] node [anchor=south] {\color{black}{$X_3$}};
		\draw[->] (8) arc [start angle = 135, end angle =-165, radius = 4mm] node [anchor=south] {\color{black}{$X_4$}};
		\draw[->] (9) arc [start angle = 135, end angle =-165, radius = 4mm] node [anchor=south] {\color{black}{$X_5$}};
		\draw[->] (10) arc [start angle = 135, end angle =-165, radius = 4mm] node [anchor=south] {\color{black}{$X_6$}};
	\end{tikzpicture}
		\caption{Cyclic coverings for $m=10$ and $r=9,r=10$}
	\end{figure}
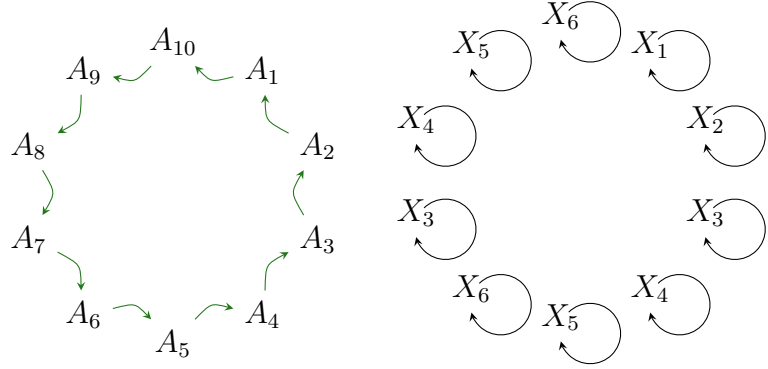

	Following the above presented color code or simply the arrows, one can trace the behavior of the $r$-cyclic operators. For example, $f$ runs through all sets of the covering, no matter where the starting point was, for $r=1$, $r=3$, $r=7$ and $r=9$. For each of $r=2$, $r=4$, $r=6$ and $r=8$ there arise two circuits. Again, the circuits are distinct, although the sets involved are the same.
	
	For example, all of the following are obviously distinct:
	\begin{align*}
		& X_1\cupdot X_3\cupdot X_5 \cupdot X_7\cupdot X_9 \text{ (when $f$ is 2-cyclic)},\\
		& X_1\cupdot X_5\cupdot X_9 \cupdot X_3\cupdot X_7 \text{ (when $f$ is 4-cyclic)},\\
		& X_1\cupdot X_7\cupdot X_3 \cupdot X_9\cupdot X_5 \text{ (when $f$ is 6-cyclic)},\\
		& X_1\cupdot X_9\cupdot X_7 \cupdot X_5\cupdot X_3 \text{ (when $f$ is 8-cyclic)}.
	\end{align*}

	The most "colorful" case is that when $f$ is 5-cyclic and the closed circuits that arise reach the number of five, each of them containing only two sets of the covering.
	
\subsection{Two lemmas on $r$-cyclic operators}\ 

As announced in the beginning of this section, our aim was to generalize what can be noticed in the behavior of $r$-cyclic operators for various values of $m$ and $r$. It is obvious that this behavior depends on the relation between the two parameters $m$ and $r$. There are two different situations, namely when $\gcd{(m,r)}=1$ and when $\gcd{(m,r)}>1$. For each of them we formulate a corresponding lemma, that will be later essential in proving our fixed point results for $r$-cyclic operators.

Let us begin with the case when $\gcd{(m,r)}=k>1$, as it generates a more spectacular behavior. Note that in the following lemmas, as well as in the fixed point results we shall state later, $r$ is no longer supposed to take also the value $m$, as the discussion would make no sense in that case.

\begin{lem} \label{Lem_mrk}
	Let $f:X\to X$ be $r$-cyclic w.r.t. $X=\overset{m}{\underset{i=1}{\bigcupdot}}X_i$, $m\geq2$, $1\leq r< m$ integers. If the greatest common divisor of $m$ and $r$, namely $\gcd(m,r)=k>1$, then there exists a covering $X=\overset{k}{\underset{i=1}{\bigcupdot}}Y_i$ with the following properties:
	\begin{itemize}
		\item[1)] The $k$ subcoverings $Y_j$, $1\leq j\leq k$ are given by:
		 			\begin{align*}
		 				& Y_1=X_1\cupdot X_{1+r} \cupdot \dots \cupdot X_{1+(\frac{m}{k}-1)r},\\	
		 				& Y_2=X_2\cupdot X_{2+r} \cupdot \dots \cupdot X_{2+(\frac{m}{k}-1)r},\\	
		 				& \vdots\\
		 				& Y_k=X_k\cupdot X_{k+r} \cupdot \dots \cupdot X_{k+(\frac{m}{k}-1)r},
		 			\end{align*}
	 			or generally $Y_j = \overset{\frac{m}{k}-1}{\underset{i=0}{\bigcupdot}}X_{j+i\cdot r}, 1\leq j \leq k$.	
	 			
	 			Note that by $X_p$ we mean $X_{p \mod m}$, each time when $p>m$.\\
	 			
	 	\item[2)] The subcoverings $Y_j$, $1\leq j\leq k$ are invariant for $f$.\\
	 	
	 	\item[3)] For each $j\in\{1,\dots,k\}$, the corresponding restriction $f|_{Y_j}$ is a cyclic operator w.r.t. $Y_j$.
	\end{itemize}
\end{lem}

\begin{remark}
	In view of Remark \ref{Rem_circuit} above, all the subcoverings $Y_j$, $1\leq j\leq k$ are (closed) circuits. 
\end{remark}

	\textbf{Comments. }So any fixed point problem for the $r$-cyclic operator $f:X\to X$ w.r.t. $\overset{m}{\underset{i=1}{\bigcupdot}}X_i$ can be "split" into $k=\gcd(m,r)$ fixed point problems for the restrictions $f|_{Y_j}:Y_j\to Y_j$, which are all cyclic operators w.r.t. $Y_j$, $1\leq j\leq k$.
	
	If $\overset{k}{\underset{j=1}{\bigcap}} Y_j=\emptyset$, then clearly $f$ will have no fixed point in $X$. 
	
	Identifying these closed circuits $Y_j$, $1\leq j\leq k$ can play an important role in certain applications. If for example one of them is causing $\overset{k}{\underset{j=1}{\bigcap}} Y_j=\emptyset$, this can be removed from the initial domain of $f$ and the fixed point problem can be studied on the remaining domain, where it might have fixed point(s). 

    The above Lemma \ref{Lem_mrk} refers to the cases when $\gcd(m,r)>1$. In the remaining case $\gcd(m,r)=1$, the associated fixed point problems are consequently less "endangered" from the point of view of the convergence of the Picard iteration, as we shall see. In practical situations, this can be either an advantage, or a disadvantage that has to be somehow avoided, depending on the nature of the cyclic problem, that is, if potential fixed points correspond to a desirable state or not.\\
    
    Having in view Lemma \ref{Lem_mrk}, the following is obvious:
    
    \begin{lem} \label{Lem_mr1}
    		Let $f:X\to X$ be $r$-cyclic w.r.t. $X=\overset{m}{\underset{i=1}{\bigcupdot}}X_i$, where $1\leq r<m$. If  $\gcd(m,r)=1$, then for any $x\in X$ the sequence $\{f^n(x)\}_{n\geq0}$ will have infinitely many terms in each $X_i$, $1\leq i\leq m$.
    \end{lem}
    
    \begin{remark} \label{Rem_FP.mr}
    	The above lemmas will play an important role when studying the fixed points of $r$-cyclic operators.
    	
    	 This is because, as mentioned at the very beginning, proving the convergence of the Picard iteration to the fixed point requires the presence of infinitely many terms of the sequence of successive approximations in each subset of the cyclic covering on which the operator is invariant.
    \end{remark}
    
	\subsection{Mathematical notation instead of visual representation. The case $m=12$}\ 
	
	Now it would be useful to verify the notations in the above lemmas on a particular case. This time we choose $m=12$ and $r=1,2,\dots,11$.  
    
    \begin{ex} \label{Ex_mare}
    	Let $f:X\to X$ be $r$-cyclic w.r.t. $X=\overset{12}{\underset{i=1}{\bigcupdot}}X_i$.
    	
    	In the sequel we shall analyze how things work for all possible values of $r$, that is, when $1\leq r<12$.
    	
    	\begin{itemize}
    		\item[1)] If $r=1$, then $f$ is cyclic w.r.t. $X=\overset{12}{\underset{i=1}{\bigcupdot}}X_i$, no further discussion is needed.
    		
    		\item[2)] If $r\in\{5,7,11\}$, then $\gcd(m,r)=1$ and no special problems are generated, in terms of Remark \ref{Rem_FP.mr}.
    		
    		\item[3)] If $r=2$, then $k=\gcd(2,12)=2$. According to Definition \ref{Def_rCycCov}, if $f$ is $2$-cyclic w.r.t. $X=\overset{12}{\underset{i=1}{\bigcupdot}}X_i$, then the following hold:
    			\begin{align*}
    			& f(X_1)\subseteq X_3, f(X_2)\subseteq X_4, f(X_3)\subseteq X_5,\dots,\\	
    			& f(x_{10})\subseteq X_{12}, f(X_{11})\subseteq X_1,f(X_{12})\subseteq X_2.	
    		\end{align*}
    		It is easy to see that there are two closed circuits, that is, two subcoverings
    			\begin{align*}
    				& Y_1=X_1\cupdot X_{3} \cupdot X_5 \cupdot X_7\cupdot X_9\cupdot X_{11},\\	
    				& Y_2=X_2\cupdot X_{4} \cupdot X_6 \cupdot X_8\cupdot X_{10}\cupdot X_{12},	
    			\end{align*}
    		which are both invariant for $f$. Moreover, according to Definition \ref{Def_CycCov}, $f|_{Y_1}:Y_1\to Y_1$ is cyclic w.r.t. $Y_1$ and $f|_{Y_2}:Y_2\to Y_2$ is cyclic w.r.t. $Y_2$.
    		
    		If we refer to the notation in Lemma \ref{Lem_mrk}, then indeed \\$1+(\frac{m}{k}-1)r=1+(\frac{12}{2}-1)\cdot 2=11$ etc.\\
    		
    		\item[4)] If $r=3$, then $k=\gcd(3,12)=3$. By Definition \ref{Def_rCycCov}, as $f$ is $3$-cyclic w.r.t. $X=\overset{12}{\underset{i=1}{\bigcupdot}}X_i$, the following hold:
    			\begin{align*}
    				& f(X_1)\subseteq X_4, f(X_2)\subseteq X_5, \dots,f(X_9)\subseteq X_{12}\\	
    				& f(x_{10})\subseteq X_{1}, f(X_{11})\subseteq X_2,f(X_{12})\subseteq X_3.	
    			\end{align*}
    		Again it is easy to see that in this case there are three subcoverings:
    			\begin{align*}
    				& Y_1=X_1\cupdot X_4 \cupdot X_7 \cupdot X_{10},\\	
    				& Y_2=X_2\cupdot X_5 \cupdot X_8 \cupdot X_{11},\\	
    				& Y_3=X_3\cupdot X_6 \cupdot X_9 \cupdot X_{12},
    			\end{align*}
    		which are all invariant for $f$ and the respective restrictions are cyclic operators w.r.t. the corresponding $Y_j$, $1\leq j\leq 3$.
    		
    		\item[5)] If $r=4$, that is, $f$ is $4$-cyclic w.r.t. $X=\overset{12}{\underset{i=1}{\bigcupdot}}X_i$, then $k=\gcd(4,12)=4$ and the subcoverings
    		\begin{align*}
    			& Y_1=X_1\cupdot X_5 \cupdot X_9,    & Y_3=X_3\cupdot X_7 \cupdot X_{11},\\	
    			& Y_2=X_2\cupdot X_6 \cupdot X_{10}, & Y_4=X_4\cupdot X_8 \cupdot X_{12}    			
    		\end{align*}
    		are invariant for $f$ and cyclic coverings w.r.t. the corresponding restrictions of $f$.
    		
    		\item[6)] If $r=6$, that is, $f$ is $6$-cyclic w.r.t. $X=\overset{12}{\underset{i=1}{\bigcupdot}}X_i$, then $k=\gcd(6,12)=6$ and the subcoverings
    		\begin{align*}
    			& Y_1=X_1\cupdot X_7, &Y_3=X_3\cupdot X_9,		&\qquad Y_5=X_5\cupdot X_{11},\\	
    			& Y_2=X_2\cupdot X_8, &Y_4=X_4\cupdot X_{10},  &\qquad Y_6=X_6\cupdot X_{12}, 			
    		\end{align*}
    		are invariant for $f$ and cyclic coverings w.r.t. the corresponding restrictions of $f$.
    		
    		\item[7)] If $r=8$, that is, $f$ is $8$-cyclic $X=\overset{12}{\underset{i=1}{\bigcupdot}}X_i$, then $k=\gcd(8,12)=4$ and things turn out more interesting, as the four invariant subcoverings are
    		\begin{align*}
    			& Y_1=X_1\cupdot X_9 \cupdot X_5,    & Y_3=X_3\cupdot X_{11} \cupdot X_{7},\\	
    			& Y_2=X_2\cupdot X_{10} \cupdot X_6, & Y_4=X_4\cupdot X_{12} \cupdot X_{8}.    			
    		\end{align*}
    		In view of Remark \ref{Rem_CupCdot} above, these subcoverings are different from those we had for $r=4$, as for example $X_1\cupdot X_5 \cupdot X_9$ is not the same as $X_1\cupdot X_9 \cupdot X_5$. 
    		
    		And again they are all cyclic w.r.t. the corresponding restrictions of $f$.
    		
    		Here it is interesting (still natural) that $f$ is cyclic on all four subcoverings $Y_j,j=\overline{1,4}$, but in each case this means it is $2$-cyclic w.r.t. the corresponding subcovering $Y_j,j=\overline{1,4}$ from the case $r=4$.
    		
    		For example, here $f$ is cyclic on $X_1\cupdot X_9\cupdot X_5$, which means it is $2$-cyclic on $X_1\cupdot X_5\cupdot X_9$.
    		
    		\item[8)] If $r=9$, that is, $f$ is $9$-cyclic w.r.t. $X=\overset{12}{\underset{i=1}{\bigcupdot}}X_i$, then $k=\gcd(9,12)=3$ and the three invariant subcoverings are 
    		\begin{align*}
    			& Y_1=X_1\cupdot X_{10} \cupdot X_{7} \cupdot X_{4},\\	
    			& Y_2=X_2\cupdot X_{11} \cupdot X_8 \cupdot X_{5},\\	
    			& Y_3=X_3\cupdot X_{12} \cupdot X_9 \cupdot X_{6}.
    		\end{align*}
    		Again we see they are different from those obtained for $r=3$. Nevertheless they are all cyclic w.r.t. the corresponding restrictions of $f$.
    		
    		We note again that $f$ is $3$-cyclic w.r.t. the subcoverings obtained for $r=3$.
    		
    		\item[9)] If finally $r=10$, the two invariant subcoverings on which $f$ is cyclic are
    		\begin{align*}
    			& Y_1=X_1\cupdot X_{11} \cupdot X_9 \cupdot X_7\cupdot X_5\cupdot X_{3},\\	
    			& Y_2=X_2\cupdot X_{12} \cupdot X_{10} \cupdot X_8\cupdot X_{6}\cupdot X_{4},	
    		\end{align*}
    		which "resemble" those for $r=2$ but are still not the same. Still $f$ is $5$-cyclic on the corresponding subcoverings from the case $r=2$.
    	\end{itemize}  
    \end{ex}

    Although unusually extended, Example \ref{Ex_mare} explains many things about how $r$-cyclic coverings work and what instruments they generate when dealing with various applications. Understanding these things simplifies the proofs of the related results and opens the way for further research.
    \bigskip

    \section{Synchronous $r$-cyclic contractions on a metric space}\label {Sect_Synchr}
    
    The first problem that anyone would think to investigate, is a possible extension of the result proved in \cite{KirkSrinVeer} and included above as Theorem \ref{Th_Kirk}, and hence an extension of the contraction principle of Banach for $r$-cyclic operators, in the framework of a metric space.
    
    Recall that in \cite{KirkSrinVeer} there were considered cyclic contractions or $1$-cyclic contractions, that is, mappings satisfying $d(F(x),F(y))\leq k\cdot d(x,y)$ for any $x\in A_i$, $y\in A_{i+1}$, where $\overset{p}{\underset{i=1}{\cup}}A_i$ was a cyclic covering w.r.t. $F$.
    
    While working on a $r$-cyclic covering, one can notice two naturally arising directions:
    \begin{itemize}
    	\item[1.] when $f$ satisfies a contraction condition that holds for any $x_i\in X_i$, $y\in X_{i+1}$, that is, for elements belonging to successive sets in the cyclic covering;
    	\item[2.] when $f$ satisfies a contraction condition that holds for any $x_i\in X_i$, $y\in X_{i+r}$, where $r$ is the same constant that defines the $r$-cyclic covering.
    \end{itemize}
    
    In order to differentiate these two classes of operators, we have chosen to call the second one \textit{synchronous} $r$-cyclic contractions, as the contraction condition has the same "rhythm" as the cyclic covering, and consequently the first class will be called \textit{asynchronous} $r$-cyclic contractions.
    
    Obviously, asynchronous $1$-cyclic contractions are the same as synchronous $1$-cyclic contractions and actually they are all cyclic contractions as studied in \cite{KirkSrinVeer}. 
    
    In spite of what one could think at a first glance, among the two classes, synchronous and asynchronous $r$-cyclic contractions, the study of the first one is the easiest to approach and the current section is dedicated to it. 
    
    We start by introducing:
    
    \begin{defn} \label{Def_rCycSyn}
    	Let $(X,d)$ be a metric space, $f:X\to X$ and  $X=\overset{m}{\underset{i=1}{\bigcupdot}}X_i$ a $r$-cyclic covering w.r.t. $f$, where $m\geq2$ and $1\leq r<m$ are integers. If there exists $c\in[0,1)$ such that for any $x\in X_i, y\in X_{i+r}, 1\leq i\leq m,$
    	\[
    	d(f(x),f(y))\leq c\cdot d(x,y), 
    	\]
    	then $f$ is called a synchronous $r$-cyclic contraction w.r.t. $\overset{m}{\underset{i=1}{\bigcupdot}}X_i$.
    \end{defn}

	If $r=1$ in the above definition, then one obtains the definition of cyclic contractions studied in \cite{KirkSrinVeer}, see Theorem \ref{Th_Kirk} above.

Now let us see which conditions must be required so that operators belonging to this class have one or more fixed points. As the above Lemmas \ref{Lem_mrk} and \ref{Lem_mr1} suggest, we shall need to analyze two different situations. A first case, when $\gcd{(m,r)}=1$, is easier to explore, as expected. The result we obtain is a generalization of Theorem \ref{Th_Kirk}, where this condition was also fulfilled, as $\gcd{(m,1)}=1$.

\begin{theorem}\label{Th_r-sync.contr.1}
	Let $(X,d)$ be a complete metric space, $m\geq2$ and $1\leq r<m$ integers, $\overset{m}{\underset{i=1}{\bigcupdot}}X_i$ a covering of $X$ with $X_i\in\mathcal{P}_{cl}(X)$, $i=\overline{1,m}$ and $c\in[0,1)$ such that $f:X\to X$ is a synchronous $r$-cyclic contraction with constant $c$ w.r.t. $\overset{m}{\underset{i=1}{\bigcupdot}}X_i$.
	
	If $\gcd{(m,r)}=1$, then $f$ has a unique fixed point in $X$, that can be obtained by means of the Picard iteration starting from any point in $X$.
\end{theorem}

\begin{proof}	

Considering the above definitions, we know from the hypothesis of the theorem that:
\begin{align*}
	&f(X_1)\subseteq X_{1+r}, f(X_2)\subseteq X_{2+r}, \dots, f(X_m)\subseteq X_{m+r}\\
	\text{and}&\\
	&d(f(x),f(y))\leq c\cdot d(x,y), \text{for any }x\in X_i, y\in X_{i+r}, 1\leq i\leq m.	
\end{align*}
We aim to show that the Picard iteration of $f$ converges to a fixed point, starting from any point in $X$. 

Therefore we consider an arbitrary point $x_0\in X$. As $X=\overset{m}{\underset{i=1}{\bigcupdot}}X_i$, there exists $l\in \{1,2,\dots,m\}$ such that $x_0\in X_l$. Because $f$ is $r$-cyclic w.r.t. $\overset{m}{\underset{i=1}{\bigcupdot}}X_i$, the terms of the Picard iteration are distributed as follows:
\[
x_1=f(x_0)\in X_{l+r}, x_2=f(x_1)\in X_{l+2r},\dots,x_n=f(x_{n-1})\in X_{l+nr},\dots,
\]
where $X_k=X_{k\mod m}$ for all $k>m$.

We have that
\[
d(x_n,x_{n+1})=d(f(x_{n-1}),f(x_n)),
\]
which, in view of the fact that $x_{n-1}\in X_{l+(n-1)r}$ and $x_n\in X_{l+nr}$, implies that
\[
d(x_n,x_{n+1})\leq c\cdot d(x_{n-1},x_n),n\geq 1.
\]
This way we obtain that, for $n\geq1$,
\[
d(x_n,x_{n+1})\leq c^n\cdot d(x_0,x_1).
\]
For $p\geq1$ we have that
\begin{align*}
	d(x_n,x_{n+p})	&\leq d(x_n,x_{n+1})+d(x_{n+1},x_{n+2})+\dots+d(x_{n+p-1},x_{n+p})\\
					&\leq c^n\cdot d(x_0,x_1)+c^{n+1}\cdot d(x_0,x_1)+\dots+c^{n+p}\cdot d(x_0,x_1)\\
					&=c^n\cdot \dfrac{1-c^p}{1-c}\cdot d(x_0,x_1),
\end{align*}
so $d(x_n,x_{n+p})\to 0$ as $n\to \infty$, which shows that $\{x_n\}_{n\geq0}$ is a Cauchy sequence in the complete metric space $(X,d)$. 

So there exists $\overline{x}\in X$ such that $\underset{n\to\infty}{\lim}x_n=\overline{x}$.

At this point it is essential to know that $\{x_n\}_{n\geq0}$ has infinitely many terms in each $X_i,i=\overline{1,m}$, and this is ensured by the condition $\gcd{(r,m)}=1$, see Lemma \ref{Lem_mr1} above. So from each $X_i$ one can extract a subsequence of $\{x_n\}_{n\geq0}$ which converges to $\overline{x}$ as well, since $(X,d)$ is complete. As $X_i$, $i=\overline{1,m}$ are all closed, it follows that 
\[
\overline{x}\in \overset{m}{\underset{i=1}{\bigcap}}X_i.
\]

Then $\overset{m}{\underset{i=1}{\bigcap}}X_i$ is not empty and is invariant for $f$. It follows that the restriction $f|_{\overset{m}{\underset{i=1}{\cap}}X_i}$ is a Banach contraction with constant $c\in[0,1)$ on the complete metric space $\overset{m}{\underset{i=1}{\bigcap}}X_i$, so it has a unique fixed point, say $x^*\in \overset{m}{\underset{i=1}{\bigcap}}X_i$, which can be obtained as the limit of the Picard iteration starting from any initial point $x\in \overset{m}{\underset{i=1}{\bigcap}}X_i.$

We still have to prove that the Picard iteration converges to $x^*$ for any starting point in $X$. So let $x\in X$. Since $X=\overset{m}{\underset{i=1}{\bigcupdot}}X_i,$ there must exist $s\in \{1,2,\dots,m\}$ such that $x\in X_s$. As $x^*\in \overset{m}{\underset{i=1}{\bigcap}}X_i$, it follows that $x^*\in X_{s+r}$ as well and, since $f$ is a synchronous $r$-cyclic contraction w.r.t. $\overset{m}{\underset{i=1}{\bigcupdot}}X_i$, we have that
\[
d(f(x),x^*)=d(f(x),f(x^*))\leq c\cdot d(x,x^*).
\]
But $f(x)\in X_{s+r}$ and $x^*=f(x^*)\in X_{s+2r}$, so 
\[
d(f^2(x),x^*)=d(f(f(x)),f(x^*))\leq c\cdot d(f(x),x^*)\leq c^2\cdot d(x,x^*).
\]
Following in a similar manner we obtain that, for $n\geq 1$,
\[
d(f^n(x),x^*)\leq c^n\cdot d(x,x^*).
\]
Then $f^n(x)\to x^*$ as $n\to\infty$, so the Picard iteration converges to the unique fixed point $x^*\in\overset{m}{\underset{i=1}{\bigcap}}X_i$ for any starting point $x\in X$.
\end{proof}

In order to complete the study of synchronous $r$-cyclic contractions, we have to analyze the remaining situations that are not covered by Theorem \ref{Th_r-sync.contr.1}, namely when $\gcd{(m,r)}>1$. Having in view Lemma \ref{Lem_mrk} and Definition \ref{Def_rCycSyn}, the following is obvious:

\begin{lem}\label{Lem_mrk.synch}
	Let $(X,d)$ be a metric space, $m\geq2$, $1\leq r<m$ and $f:X\to X$ a synchronous $r$-cyclic contraction w.r.t. $X=\overset{m}{\underset{i=1}{\bigcupdot}}X_i$. 
	
	If $\gcd{(m,r)}=k>1$, then there exists a covering $X=\overset{k}{\underset{i=1}{\bigcupdot}}Y_i$ with the following properties:
	\begin{itemize}
		\item[1)] 
		$Y_j = \overset{\frac{m}{k}-1}{\underset{i=0}{\bigcupdot}}X_{j+i\cdot r}, 1\leq j \leq k$.	\\
		
		\item[2)] The subcoverings $Y_j$, $1\leq j\leq k$ are invariant for $f$.\\
		
		\item[3)] For each $j\in\{1,\dots,k\}$, the corresponding restriction of $f$ is a cyclic contraction (or $1$-cyclic contraction) w.r.t. $Y_j$.
	\end{itemize} 
\end{lem}    

Now based on the above results we have:

\begin{theorem}\label{Th_r-sync.contr.k}
	Let $(X,d)$ be a complete metric space, $m\geq2$ and $1\leq r<m$ integers, $\overset{m}{\underset{i=1}{\bigcupdot}}X_i$ a covering of $X$ with $X_i\in\mathcal{P}_{cl}(X)$, $i=\overline{1,m}$ and $c\in[0,1)$ such that $f:X\to X$ is synchronous $r$-cyclic contraction with constant $c$ w.r.t. $\overset{m}{\underset{i=1}{\bigcupdot}}X_i$.
	
	If $\gcd{(m,r)}=k>1$, then $f$ is a weakly Picard operator and it has at most $k$ fixed points $x^*_j\in Y_j$, where $Y_j = \overset{\frac{m}{k}-1}{\underset{i=0}{\bigcupdot}}X_{j+i\cdot r}, 1\leq j \leq k$.
	
	For each $1\leq j\leq k$, $x^*_j$ can be obtained as the limit of the Picard iteration of $f$, starting from any initial point in $Y_j$.	
\end{theorem}

\begin{proof}
	By Lemma \ref{Lem_mrk.synch} and applying the Theorem \ref{Th_Kirk} due to Kirk et al. \cite{KirkSrinVeer} for the restrictions $f|_{Y_j},j=\overline{1,k}$.
\end{proof}

	As we see in Theorem \ref {Th_r-sync.contr.1}, the convergence of the Picard iteration to the unique fixed point of a synchronous $r$-cyclic contraction w.r.t. a covering $\overset{m}{\underset{i=1}{\bigcupdot}}X_i$ of a complete metric space $X$ is guaranteed only if $\gcd{(m,r)}=1$.
	
	In the other case, that is, when $\gcd{(m,r)}=k>1$, the fixed point problem can be approached by splitting it into $k$ fixed point problems for the restrictions $f|_{Y_j}, j=\overline{1,k}$, see Lemma \ref{Lem_mrk} and Theorem \ref{Th_r-sync.contr.k}. These restrictions of $f$ are $1$-cyclic or simply cyclic contractions with the same constant $c\in[0,1)$ on the respective subcoverings of "length" $\frac{m}{k}$, so there are known instruments available to study the existence of the fixed points there.

\bigskip

\section{Some remarks on synchronous $r$-cyclic contractions}

\bigskip
Many questions may arise now regarding the above results. We shall analyze them under the assumptions of Theorem \ref{Th_r-sync.contr.k} above.

What happens if the subcoverings $Y_j, j=\overline{1,k}$ are pairwise disjoint? What happens if two or many of them have at least a common element, even if $\overset{k}{\underset{j=1}{\bigcap}}Y_j=\emptyset$? If $f$ acts independently on the $k$ circuits $Y_j, j=\overline{1,k}$, then which consequences might this have? 

In the following we try to give some answers.

\begin{remark}\label{Rem_Yij.1}
	First, if the subcoverings $Y_j,j=\overline{1,k}$ are pairwise disjoint, that is, $Y_i\cap Y_j=\emptyset$, for $i\neq j$, $i,j=\overline{1,k}$, then certainly $f$ will have the maximum number of fixed points, namely $k$ fixed points, one in each $Y_j$, $j=\overline{1,k}$. This stands for the case when the operator $f$ splits its domain in $k$ separate circuits $Y_j=\overset{\frac{m}{k}-1}{\underset{i=0}{\bigcup}}X_{j+ir}$, being 1-cyclic on each of them.
\end{remark}

\begin{remark}	
	In view of Definition 4.4 from \cite{Rus-Set-Th.CJM}, under the assumptions of Theorem \ref{Th_r-sync.contr.k} above, the covering $\overset{k}{\underset{j=1}{\bigcup}}Y_j$ is a fixed point invariant partition of $X$.
\end{remark}

\begin{remark}\label{Rem_Yij.2}
	Now let us see what happens if $Y_j$, $j=\overline{1,k}$ are not pairwise disjoint, that is, there are at least two having a common element. We may assume without loss of generality that $Y_1\cap Y_2\neq\emptyset$, so there is some $\overline{x}\in Y_1\cap Y_2$.
	
	Again without loss of generality, due to the properties of a cyclic covering, we may assume that $\overline{x}\in X_1$ and $\overline{x}\in X_2$ as well, since 
	 \begin{align*}
	 	& Y_1=X_1\cupdot X_{1+r} \cupdot \dots \cupdot X_{1+(\frac{m}{k}-1)r},\\	
	 	& Y_2=X_2\cupdot X_{2+r} \cupdot \dots \cupdot X_{2+(\frac{m}{k}-1)r}.
	 \end{align*}
  
  	Then $f(\overline{x})\in X_{1+r}\cap X_{2+r}$, $f^2(\overline{x})\in X_{1+2r}\cap X_{2+2r}$ and so on.
  	
  	But since $f$ is $r$-cyclic on $X$, it follows that the Picard iteration starting from $\overline{x}$ will converge to the same fixed point $\overline{x}^*$ of $f$ which belongs to $\overset{\frac{m}{k}-1}{\underset{i=0}{\bigcap}}X_{1+ir}$ and to $\overset{\frac{m}{k}-1}{\underset{i=0}{\bigcap}}X_{2+ir}$ as well.
  	
  	Practically this says that if we determine that two such circuits on which $f$ is cyclic have at least one common point, then they have a lot of common points and $f$ will have a unique fixed point in $Y_1\cap Y_2$, which can be obtained as the limit of the Picard iteration starting from any point in $Y_1\cup Y_2$, even if this initial point is not a common point of $Y_1$ and $Y_2$.
  	
  	This remains valid if we talk about more than two circuits for which one can identify a common element.
\end{remark}

\begin{remark}\label{Rem_Yij.3}
	Consequently, if there is at least one element common to all $Y_j$, $1\leq j\leq m$, that is, $\overset{k}{\underset{j=1}{\bigcap}}Y_j\neq\emptyset$, then $f$ will have only one fixed point in $\overset{m}{\underset{i=1}{\bigcap}}X_i\neq\emptyset$, which can be obtained as the limit of the Picard iteration starting from any point in $X$.
\end{remark}

\begin{remark}
	If $f$ is 1-cyclic on $\overset{m}{\underset{i=1}{\bigcupdot}}X_i$ with $\overset{m}{\underset{i=1}{\bigcap}}X_i=\emptyset$, we know for sure it has no fixed point.

	If $X=\overset{m}{\underset{i=1}{\bigcupdot}}X_i$ is $r$-cyclic w.r.t. $f:X\to X$ and $\gcd{(m,r)}=k>1$, then the fact that $\overset{k}{\underset{j=1}{\bigcap}}Y_j=\emptyset$ (see Lemma \ref{Lem_mrk} above) necessarily implies that $\overset{m}{\underset{i=1}{\bigcap}}X_i=\emptyset$ (but not that $\overset{m}{\underset{i=1}{\bigcup}}X_i$ is a partition of $X$), while the converse doesn't have to hold generally.
	
	Indeed, when $(X,d)$ is a complete metric space, $X_i\in\mathcal{P}_{cl}(X)$ and $f$ is a synchronous $r$-cyclic contraction on $X=\overset{m}{\underset{i=1}{\bigcupdot}}X_i$ with $\overset{m}{\underset{i=1}{\bigcap}}X_i=\emptyset$ and $r\geq2$, then $f$ could have several fixed points, one in each subcovering $Y_j$ with nonempty intersection $\overset{\frac{m}{k}-1}{\underset{i=0}{\bigcap}}X_{j+i\cdot r}$.
\end{remark}

\bigskip

	
 \section{Asynchronous $r$-cyclic contractions}

As we have already mentioned above, the cyclic contraction condition studied in \cite{KirkSrinVeer} can be extended in two directions if working on $r$-cyclic coverings with $r\geq2$. The first direction was explored in Section \ref{Sect_Synchr}, by defining and studying the synchronous $r$-cyclic contractions, namely those cyclic operators which satisfy a contraction condition for $x\in X_i$, $y\in X_{i+r}$, $i=\overline{1,m}$.

In the current section we will explore the case when the contraction condition holds for $x\in X_i$, $y\in X_{i+1}$, $i=\overline{1,m}$, that is, for elements belonging to successive sets in the cyclic covering. Therefore we introduce:

    \begin{defn} \label{Def_rCycCons}
	Let $(X,d)$ be a metric space, $f:X\to X$ and  $X=\overset{m}{\underset{i=1}{\bigcupdot}}X_i$ a $r$-cyclic covering w.r.t. $f$, where $m\geq2$ and $1\leq r<m$ are integers. If there exists $c\in[0,1)$ such that for any $x\in X_i, y\in X_{i+1}, i=\overline{1,m},$
	\begin{equation*}\label{Ec_rCycCons}
	d(f(x),f(y))\leq c\cdot d(x,y), 
	\end{equation*}
	then $f$ is called an asynchronous $r$-cyclic contraction w.r.t. $\overset{m}{\underset{i=1}{\bigcupdot}}X_i$.
\end{defn}

\begin{remark}	
	Because for $r=1$ the discussion would reduce to the cyclic contraction condition in \cite{KirkSrinVeer}, we shall continue the discussion for $r\geq2$. Consequently we shall take $m\geq3$.
\end{remark}

The behavior of the asynchronous $r$-cyclic contractions is not more complicated than that of the synchronous $r$-cyclic contractions and besides for both cases the things work pretty intuitive. But if one wishes to keep the notation rigorous for the general case, this becomes a little complicated.

That is why we shall first analyze the case when $r=2$. As in the case of synchronous $r$-cyclic contractions, where we had two separate fixed point results, the discussion has to cover two situations: when $\gcd{(m,r)}=1$ and when $\gcd{(m,r)}=k>1$. 

For the case $r=2$ this means $m$ odd or even, and both situations are covered in the next result:

\begin{theorem}\label{Th_r-cons.contr}
	Let $(X,d)$ be a complete metric space, $X=\overset{m}{\underset{i=1}{\bigcupdot}}X_i$, $m\geq3$ a closed covering and $c\in[0,1)$ such that $f:X\to X$ is an asynchronous $2$-cyclic contraction with constant $c$ w.r.t. $\overset{m}{\underset{i=1}{\bigcupdot}}X_i$.
	
	Then $f$ is a Picard operator.
\end{theorem}

\begin{proof}
	
	Let us start by resuming what the hypothesis of the theorem actually implies. Since $f$ is 2-cyclic operator on $\overset{m}{\underset{i=1}{\cupdot}}X_i$, we have that
		\[
			f(X_1)\subseteq X_3, f(X_2)\subseteq X_4, \dots, f(X_{m-2})\subseteq X_m, f(X_{m-1})\subseteq X_1,f(X_m)\subseteq X_2.
		\]
	The fact that $f$ is asynchronous cyclic contraction means that
		\[
			d(f(x),f(y))\leq c\cdot d(x,y),
		\]
	for any $x\in X_i,y\in X_{i+1},i=\overline{1,m}$.
	
	In this case we need to begin our proof by considering two orbits starting from two different points $x_0$ and $y_0$ belonging to successive sets of the cyclic covering. Without loss of generality, we may assume that $x_0\in X_1$ and $y_0\in X_2$.
	
	The orbits $\{f^n(x_0)\}_{n\geq0}$ and $\{f^n(y_0)\}_{n\geq0}$ have their points in $X_1, X_3,X_5,\dots$ and, respectively, in $X_2, X_4, X_6, \dots$.
	
	It follows that
	\begin{align*}
		& d(x_1,y_1)=d(f(x_0),f(y_0))\leq c\cdot d(x_0,y_0)\\
		& d(y_1,x_2)=d(f(y_0),f(x_1))\leq c\cdot d(y_0,x_1),
	\end{align*}
	then
	\begin{align*}
		& d(x_2,y_2)=d(f(x_1),f(y_1))\leq c^2\cdot d(x_0,y_0)\\
		& d(y_2,x_3)=d(f(y_1),f(x_2))\leq c^2\cdot d(y_0,x_1),
	\end{align*}
	and so on
	\begin{align*}
		 &d(x_n,y_n) \leq c^n\cdot d(x_0,y_0)\\
		 &d(y_n,x_{n+1})\leq c^n\cdot d(y_0,x_1).\nonumber		 
	\end{align*}
	So for $n\geq1$ we have that
	\[
	d(x_n,x_{n+1})\leq d(x_n,y_n)+d(y_n,x_{n+1})	\leq c^n\cdot[d(x_0,y_0)+d(y_0,x_1)].
	\]
	If we denote by $A=d(x_0,y_0)+d(y_0,x_1)\geq0$, we have that
	\[
	d(x_n,x_{n+1})\leq c^n\cdot A, n\geq1.
	\]
	For $p\geq1$ we obtain that
	\[
	d(x_n,x_{n+p})\leq c^n\cdot\dfrac{1-c^p}{1-c}\cdot A,
	\]
	which leads to the conclusion that $\{x_n\}_{n\geq0}$ is a Cauchy sequence in the complete metric space $(X,d)$, so there exists its limit $\overline{x}\in X$.
	
	At this point we have to analyze the two different cases mentioned above.
	
	\begin{itemize}
		\item[1)] If $m$ is odd, so $\gcd{(m,2)}=1$, the sequence $\{x_n\}_{n\geq0}$ has infinitely many terms in each $X_i$, $i=\overline{1,m}$ (see Lemma \ref{Lem_mr1}), so from each $X_i$ one can extract a subsequence of $\{x_n\}_{n\geq0}$ that converges to $\overline{x}$. 
		
		As $X_i$, $i=\overline{1,m}$ are all closed, it follows that $\overline{x}\in\overset{m}{\underset{i=1}{\bigcap}}X_i$.

		\item[2)] If $m$ is even, so $\gcd{(m,2)}=2$, there are two circuits
			 \[X_1\cupdot X_3\cupdot \dots\cupdot X_{m-1} \text{ and } X_2\cupdot X_4\cupdot\dots\cupdot X_{m} \]
		which are invariant for $f$ (see Lemma \ref{Lem_mrk}). The problem here is that the two orbits $\{f^n(x_0)\}_{n\geq0}$ and $\{f^n(y_0)\}_{n\geq0}$ considered above do not have infinitely many terms in each $X_i$, $i=\overline{1,m}$. The first will converge in $\overset{\frac{m}{2}}{\underset{i=1}{\bigcup}}X_{2i-1}$ to $\overline{x}\in\overset{\frac{m}{2}}{\underset{i=1}{\bigcap}}X_{2i-1}$ and the second one will converge in $\overset{\frac{m}{2}}{\underset{i=1}{\bigcup}}X_{2i}$ to $\overline{y}\in\overset{\frac{m}{2}}{\underset{i=1}{\bigcap}}X_{2i}$.
		
		Still the inequality $d(x_n,y_n) \leq c^n\cdot d(x_0,y_0),n\geq0$, obtained earlier in the proof, ensures that the two orbits $\{f^n(x_0)\}_{n\geq0}$ and $\{f^n(y_0)\}_{n\geq0}$ have the same limit $\overline{x}=\overline{y}\in\overset{m}{\underset{i=1}{\bigcap}}X_i$.
	\end{itemize}

	So for both cases $m$ odd and $m$ even we know for sure that $\overset{m}{\underset{i=1}{\bigcap}}X_i\neq\emptyset$. Then the restriction $f|_{\overset{m}{\underset{i=1}{\cap}}X_i}$ is a Banach contraction with constant $c$ on $\overset{m}{\underset{i=1}{\bigcap}}X_i$, having a unique fixed point $x^*\in\overset{m}{\underset{i=1}{\bigcap}}X_i$ which can be obtained by means of the Picard iteration starting from any point in $\overset{m}{\underset{i=1}{\bigcap}}X_i$.
	
	We still have to see if $x^*$ can be obtained for any initial point in $X$. Therefore we take an arbitrary $x\in X$. Then there is some $l\in\{1,2,\dots,m\}$ such that $x\in X_l$. As $x^*\in X_{l+1}$, the following holds:
	\[
	d(f(x),f(x^*))\leq c\cdot d(x,x^*).
	\]
	Since $f$ is $2$-cyclic, it follows that $f(x)\in X_{l+2}$. But $x^*\in X_{l+3}$ as well, so
	\[
	d(f^2(x),x^*)=d(f(f(x)),f(x^*))\leq c\cdot d(f(x),x^*)\leq c^2\cdot d(x,x^*).
	\]
	In this manner we obtain that
	\[
	d(f^n(x),x^*)\leq c^n\cdot d(x,x^*), n\geq1.
	\]
	Now this implies that $f^n(x)\to x^*$, as $n\to\infty$, for any $x\in X$.
	
	In conclusion $f$ has a unique fixed point, no matter if $m$ is odd or even and this fixed point may be obtained as the limit of the Picard iteration starting from any point in $X=\overset{m}{\underset{i=1}{\bigcup}}X_i$.

\end{proof}

\begin{remark}
	Analyzing the proof of Theorem \ref{Th_r-cons.contr} it is not difficult to imagine how the same pattern works for any $r\geq2$. 
	
	Some intuitive representations are given below for $r=2$ (the two orbits considered in the proof may be seen in the figure), $r=3$ (where three orbits starting from three consecutive sets of the covering have to be considered) and $r=4$ (with four such orbits). 
	
	In each case the proof would be similar to the one of Theorem \ref{Th_r-cons.contr}.\\

	\begin{itemize}
		\item[1)] The case $r=2$, as in Theorem \ref{Th_r-cons.contr}. As one can see, both orbits $\{x_n\}_{n\geq0}$ and $\{y_n\}_{n\geq0}$ converge to the same point, with: 
		\begin{align*}
			& \{x_n\}_{n\geq0}\subseteq X_1\cup X_3\cup X_5 \cup \dots\\
			& \{y_n\}_{n\geq0}\subseteq X_2\cup X_4\cup X_6 \cup \dots,
		\end{align*}
		
		while $d(x_n,y_n)\to0\text{ as }n\to \infty.$\\
		\begin{figure}[H]
		\begin{center}
			\begin{tikzpicture}[scale=7]
			\def\a{8/9};
			\coordinate[label={[myblue]left:$x_0$}] (A) at (0,0);
			\coordinate[label={[plum]right:$y_0$}] (B) at (1,0);
			
			\draw[gray!80,-] (A) to (B); 
			
			\coordinate[label] (X) at (A); 	\coordinate[label] (Y) at (B);
			
			\foreach \i in {1,...,4}
			{
				\coordinate[label={[myblue]left: \small{$x_{\i}$}}] (C) at ($ (Y)!\a!7:(X) $);
				\draw[gray!80,-]  (Y) to (C);
				
				\coordinate[label={[plum]right:\small{$y_{\i}$}}] (D) at ($ (C)!\a!-8:(Y) $);
				\draw[gray!80,-] (C) to (D);
				
				\coordinate[label] (X) at (C);
				\coordinate[label] (Y) at (D);
			}
			\foreach \i in {5,...,8}
			{
				\coordinate[label={[myblue]left: \tiny{$x_{\i}$}}] (C) at ($ (Y)!\a!10:(X) $);
				\draw[gray!80,-]  (Y) to (C);
				
				\coordinate[label={[plum]right:\tiny{$y_{\i}$}}] (D) at ($ (C)!\a!-10:(Y) $);
				\draw[gray!80,-] (C) to (D);
				
				\coordinate[label] (X) at (C);
				\coordinate[label] (Y) at (D);
			}
			\foreach \i in {1,...,10}
			{
				\coordinate (C) at ($ (Y)!\a!10:(X) $);
				\draw[gray!80,-]  (Y) to (C);
				
				\coordinate (D) at ($ (C)!\a!-10:(Y) $);
				\draw[gray!80,-] (C) to (D);
				
				\coordinate[label] (X) at (C);
				\coordinate[label] (Y) at (D);
			}
		\end{tikzpicture} 
	\caption{Convergence of the two sequences in the case of a 2-cyclic asynchronous contraction}
	\end{center}
	\end{figure}
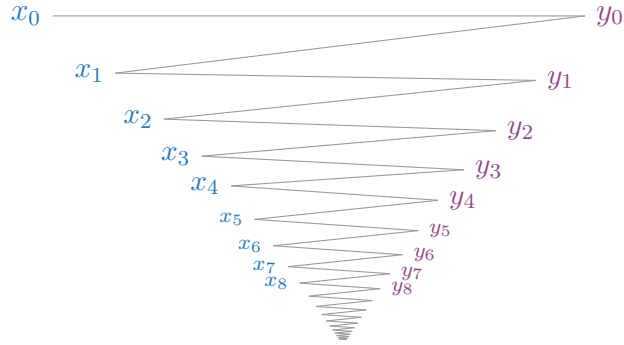

	\ \\
	\item[2)] The case $r=3$. Then the three orbits $\{x_n\}_{n\geq0}$, $\{y_n\}_{n\geq0}$ and $\{z_n\}_{n\geq0}$ converge to the same point, with: 
		\begin{align*}
			& \{x_n\}_{n\geq0}\subseteq X_1\cup X_4\cup X_7 \cup \dots \\
			& \{y_n\}_{n\geq0}\subseteq X_2\cup X_5\cup X_8 \cup \dots \\
			& \{z_n\}_{n\geq0}\subseteq X_3\cup X_6\cup X_9 \cup \dots,
		\end{align*}
	while
		\begin{align*}
			& d(x_n,y_n)\to 0\text{ as }n\to \infty\\
			& d(y_n,z_n)\to 0\text{ as }n\to \infty
		\end{align*}
	and consequently
		\[d(x_n,z_n)\to 0\text{ as }n\to \infty.\\
		\]
	This can be visualized in Figure \ref{Fig_3-asynch}.
	\begin{figure}
		\begin{center}
			\begin{tikzpicture}[scale=6]
			\def\a{8/9};
			
			\coordinate[label={[myblue]left:$x_0$}] (A) at (225-360/3:1);
			\coordinate[label={[plum]right:$y_0$}] (B) at (225-2*360/3:1);
			\coordinate[label={[mygreen]below:$z_0$}] (C) at (225-3*360/3:1);
			
			\draw[gray!80,-] (A) to (B); \draw[gray!80,-] (B) to (C);
			
			\coordinate (X) at (A); 	\coordinate (Y) at (B); \coordinate (Z) at (C);
			
			\foreach \i in {1,...,4}
			{
				\coordinate[label={[myblue]right: \small{$x_{\i}$}}] (P) at ($ (Z)!\a!-6:(X) $);
				\draw[gray!80,-]  (Z) to (P);
				
				\coordinate[label={[plum]below:\small{$y_{\i}$}}] (Q) at ($ (X)!\a!-8:(Y) $);
				\draw[gray!80,-] (P) to (Q);
				
				\coordinate[label={[mygreen]left:\small{$z_{\i}$}}] (R) at ($ (Y)!\a!-8:(Z) $);
				\draw[gray!80,-] (Q) to (R);
				
				\coordinate (X) at (P);
				\coordinate (Y) at (Q);
				\coordinate (Z) at (R);
			}
			\foreach \i in {5,...,8}
			{
				\coordinate[label={[myblue]right: \tiny{$x_{\i}$}}] (P) at ($ (Z)!\a!-2:(X) $);
				\draw[gray!80,-]  (Z) to (P);
				
				\coordinate[label={[plum]below:\tiny{$y_{\i}$}}] (Q) at ($ (X)!\a!-3:(Y) $);
				\draw[gray!80,-] (P) to (Q);
				
				\coordinate[label={[mygreen]left:\tiny{$z_{\i}$}}] (R) at ($ (Y)!\a!-5:(Z) $);
				\draw[gray!80,-] (Q) to (R);
				
				\coordinate (X) at (P);
				\coordinate (Y) at (Q);
				\coordinate (Z) at (R);
			}
			\foreach \i in {1,...,6}
			{
				\coordinate(P) at ($ (Z)!6/8!-2:(X) $);
				\draw[gray!80,-]  (Z) to (P);
				
				\coordinate (Q) at ($ (X)!\a!-3:(Y) $);
				\draw[gray!80,-] (P) to (Q);
				
				\coordinate (R) at ($ (Y)!\a!-5:(Z) $);
				\draw[gray!80,-] (Q) to (R);
				
				\coordinate (X) at (P);
				\coordinate (Y) at (Q);
				\coordinate (Z) at (R);
			}
		\end{tikzpicture}
		\caption{Convergence of the three sequences in the case of a 3-cyclic asynchronous contraction}\label{Fig_3-asynch}
		\end{center}
	\end{figure}

	\item[3)] The case $r=4$. The four orbits $\{x_n\}_{n\geq0}$, $\{y_n\}_{n\geq0}$, $\{z_n\}_{n\geq0}$ and $\{t_n\}_{n\geq0}$ converge to the same point, while the following hold: 
	\begin{align*}
		& \{x_n\}_{n\geq0}\subseteq X_1\cup X_5\cup X_9 \cup \dots \\
		& \{y_n\}_{n\geq0}\subseteq X_2\cup X_6\cup X_{10} \cup \dots \\
		& \{z_n\}_{n\geq0}\subseteq X_3\cup X_7\cup X_{11} \cup \dots \\
		& \{t_n\}_{n\geq0}\subseteq X_4\cup X_8\cup X_{12} \cup \dots \\
	\end{align*}
while
	\begin{align*}
		& d(x_n,y_n)\to 0\text{ as }n\to \infty\\
		& d(y_n,z_n)\to 0\text{ as }n\to \infty\\
		& d(z_n,t_n)\to 0\text{ as }n\to \infty
	\end{align*}
and consequently
	\[
		 d(t_n,x_n)\to 0\text{ as }n\to \infty.
	\]
This can be visualized in Figure \ref{Fig_4-asynch}.
	\begin{figure}
	\begin{center}
			\begin{tikzpicture}[scale=6]
			\def\a{11/13};
			
			\coordinate[label={[myblue]left:$x_0$}] (A) at (180-360/4:1);
			\coordinate[label={[plum]right:$y_0$}] (B) at (180-2*360/4:1);
			\coordinate[label={[mygreen]below:$z_0$}] (C) at (180-3*360/4:1);
			\coordinate[label={[tangelo]below:$t_0$}] (D) at (180-4*360/4:1);
			
			\draw[gray!80,-] (A) to (B); \draw[gray!80,-] (B) to (C); \draw[gray!80,-] (C) to (D);
			
			\coordinate (X) at (A);
			\coordinate (Y) at (B); 
			\coordinate (Z) at (C);
			\coordinate (T) at (D);
			
			\foreach \i in {1,...,4}
			{
				\coordinate[label={[myblue]left: \small{$x_{\i}$}}] (P) at ($ (T)!\a!-5:(X) $);
				\draw[gray!80,-]  (T) to (P);
				
				\coordinate[label={[plum]right:\small{$y_{\i}$}}] (Q) at ($ (X)!\a!-5:(Y) $);
				\draw[gray!80,-] (P) to (Q);
				
				\coordinate[label={[mygreen]right:\small{$z_{\i}$}}] (R) at ($ (Y)!\a!-5:(Z) $);
				\draw[gray!80,-] (Q) to (R);
				
				\coordinate[label={[tangelo]right:\small{$t_{\i}$}}] (S) at ($ (Z)!\a!-5:(T) $);
				\draw[gray!80,-] (R) to (S);
				
				\coordinate (X) at (P);
				\coordinate (Y) at (Q);
				\coordinate (Z) at (R);
				\coordinate (T) at (S);
			}
			\def\a{7/9};
			\foreach \i in {5,...,10}
			{
				\coordinate[label={[myblue]left: \tiny{$x_{\i}$}}] (P) at ($ (T)!\a!-5:(X) $);
				\draw[gray!80,-]  (T) to (P);
				
				\coordinate[label={[plum]right:\tiny{$y_{\i}$}}] (Q) at ($ (X)!\a!-5:(Y) $);
				\draw[gray!80,-] (P) to (Q);
				
				\coordinate[label={[mygreen]right:\tiny{$z_{\i}$}}] (R) at ($ (Y)!\a!-5:(Z) $);
				\draw[gray!80,-] (Q) to (R);
				
				\coordinate[label={[tangelo]right:\tiny{$t_{\i}$}}] (S) at ($ (Z)!\a!-5:(T) $);
				\draw[gray!80,-] (R) to (S);
				
				\coordinate (X) at (P);
				\coordinate (Y) at (Q);
				\coordinate (Z) at (R);
				\coordinate (T) at (S);
			}
			
			\foreach \i in {1,...,8}
			{
				\coordinate (P) at ($ (T)!\a!-5:(X) $);
				\draw[gray!80,-]  (T) to (P);
				
				\coordinate (Q) at ($ (X)!\a!-5:(Y) $);
				\draw[gray!80,-] (P) to (Q);
				
				\coordinate (R) at ($ (Y)!\a!-5:(Z) $);
				\draw[gray!80,-] (Q) to (R);
				
				\coordinate (S) at ($ (Z)!\a!-5:(T) $);
				\draw[gray!80,-] (R) to (S);
				
				\coordinate (X) at (P);
				\coordinate (Y) at (Q);
				\coordinate (Z) at (R);
				\coordinate (T) at (S);
			}
			
		\end{tikzpicture}
		\caption{Convergence of the four sequences in the case of a 4-cyclic asynchronous contraction}\label{Fig_4-asynch}
		\end{center}
	\end{figure}
	\end{itemize}
\end{remark}

    \ \\Having in view the above considerations, we may state now the following theorem.
    
\begin{theorem}
	Let $(X,d)$ be a complete metric space, $X=\overset{m}{\underset{i=1}{\bigcupdot}}X_i$, $m\geq3$ a closed covering of $X$. If $f:X\to X$ is an asynchronous $r$-cyclic contraction  w.r.t. $\overset{m}{\underset{i=1}{\bigcupdot}}X_i$, $1\leq r<m$, then $f$ has a unique fixed point in $X$ that can be obtained by means of the Picard iteration starting from any point in $X$.
\end{theorem}  

\begin{proof}
	We include here only a sketch of the proof and omit the extended one, which would imply a much too complicated notation for the sake of rigor. 	
	It follows the idea in the proof of Theorem \ref{Th_r-cons.contr}, that is, one has to consider $r$ orbits starting from $r$ consecutive sets in the covering, let us say $X_1$, $X_2$, \dots, $X_r$. Using Lemmas \ref{Lem_mrk} and \ref{Lem_mr1} it is shown that $\overset{m}{\underset{i=1}{\bigcap}}X_i\neq\emptyset$, for $\gcd{(m,r)}\geq1$. The rest of the proof follows similarly. Each of the $r$ orbits converges to the same unique fixed point of $f$. 
\end{proof}  

\begin{remark}
	In the case $r=1$, the above theorem reduces to  Theorem \ref{Th_Kirk} due to Kirk et al. \cite{KirkSrinVeer}. But what does it offer in the case  $r>1$? 
	
	Practically, if  $r>1$ and $f$ is asynchronous $r$-cyclic w.r.t. $\overset{m}{\underset{i=1}{\bigcupdot}}X_i$, one may obtain the unique fixed point by making $r$ times less iterations, since any orbit $\{f^n(x)\}_{n\geq0}$ will "jump" over $r-1$ usual iterations towards the fixed point.
	
	If our fixed point problem is attached to a cyclic system having $r$ states, and there can be observed a sort of contraction at transfer from any state to the next one, then the solution (or equilibrium) can be obtained by iteration through only one of these $r$ states, as all of them would lead in a similar number of iterations to the same solution. One should simply choose out of the $r$ states the one for which the calculations or measurements are more convenient.
\end{remark}

\section{Instead of a conclusion. The cyclic operators that have been missing}

\medskip

There are still many things to inquire about $r$-cyclic operators. The first question one would probably still want to ask is if they are really necessary and not a trivial generalization. Instead of a long discourse, we could answer by the following:

\begin{prop}
	Let $X$ be a nonempty set and $X=\overset{m}{\underset{i=1}{\bigcupdot}}X_i$, $m\geq 2$ a covering. We denote by $CycG_m$ the set of $r$-cyclic operators w.r.t. $\overset{m}{\underset{i=1}{\bigcupdot}}X_i$, $1\leq r\leq m$:
	
	\[
	CycG_m = \{f_r \text{ is } r\text{-cyclic w.r.t. }\overset{m}{\underset{i=1}{\cupdot}}X_i | 1\leq r\leq m\}.
	\]
	
	Then $(CycG_m,\circ)$ is an abelian group, where $\circ$ is the composition of functions.
\end{prop}

\begin{proof}
	For any $m\geq 2$, the identity element of the $CycG_m$ group is $f_m$, which maps each $X_i$ into itself, $i=\overline{1,m}$.
	
	The associativity is easy to check with the definition.
	
	For each $f_r$, $r=\overline{1,m}$, the inverse element is $f_{m-r}$.
	
	Obviously the commutativity is also satisfied.
\end{proof}

This simple observation enables us to state that $r$-cyclic operators naturally complete the framework in which the study of cyclic operators should be carried out and that they have been somehow missing from all the previous research regarding fixed points or best proximity points for various types of cyclic operators.

\medskip
When it comes to applications, there are countless research areas which involve cyclic processes - in fact some of these have inspired the present study. 

It would be interesting to analyze some typical problems occurring in each of these areas and to see in which of them the cyclic phenomena are likely to be approached by the results given above. At a common search in the Scopus database, there are almost 700 000 papers matching the word \textit{cyclic} in their title, abstract or keywords, belonging to the most diverse areas, like engineering, biochemistry, genetics and molecular biology, chemistry and chemical engineering, neuroscience, computer science, ecosystems science, system engineering, sociology, even linguistics. Of course, not all of the phenomena referred there are matching the mathematics described in the present paper, just the term "cyclic" is not all of it, but among them there are still a significant number where our results could serve as instruments for some steps forward in the related research, because they reflect some natural behavior that can be easily observed in real life.

Another aspect worth mentioning is that our attention should not be attracted exclusively by quantitative applications, where numerical instruments work best, but also by applications where conceptual models work as well, as mathematics is a language and it consequently can well serve for formulating such models.

In what our results are concerned, they would fit in practical situations where some cyclic multi-phase processes can be identified, and where the equilibrium states are connected to the fixed points of some operators that show a contractive behavior of synchronous or asynchronous type. The fixed point or equilibrium can be attained by iteration, and the interesting fact is that under certain assumptions this can happen by considering values calculated for a single phase out of all phases, in those cases where all orbits would converge to the same fixed point. Besides, according to our results, the existence of one ore more solutions could be controlled by actions on the sets of values corresponding to each phase. 

There arise a lot of open problems: could the same type of investigation lead to relevant results if other types of contraction conditions, generalized metrics are considered, if best proximity points are studied etc.? This remains for future study, along with all the other aspects regarding $r$-cyclic operators, that have not been covered in the present paper.

\end{document}